\documentclass[a4paper, 12pt]{article}

\usepackage{amssymb,amsbsy,amsmath,amsfonts,amssymb,amscd}
\usepackage[francais]{babel}
\usepackage[T1]{fontenc}
\usepackage[latin1]{inputenc}

\newtheorem{theorem}{Th\'eor\`eme}[section]
\newtheorem{lemma}[theorem]{Lemme}
\newtheorem{proposition}[theorem]{Proposition}
\newtheorem{corollary}[theorem]{Corollaire}

\newtheorem{remark}[theorem]{Remarque}

\newtheorem{hypothese}[theorem]{Hypoth\`ese}

\newtheorem{notation}[theorem]{Notations}

\newcommand{\ds}{\displaystyle}
\newcommand{\ofr}{{\mathfrak o}}
\newcommand{\pfr}{{\mathfrak p}}

\newcommand{\lra}{\longrightarrow}
\newcommand{\noi}{\noindent}

\newcommand{\St}{{\mathbf S}{\mathbf t}}

\def\cind#1#2{\hbox{\rm c-ind}_{#1}^{#2}}
\newcommand{\ZZ}{\mathbb Z}
\newcommand{\CC}{\mathbb C}
\newcommand{\GG}{\mathbb G}

\newcommand{\HH}{\mathcal H}
\newcommand{\VV}{\mathcal V}
\newcommand{\WW}{\mathcal W}

\newcommand{\Afr}{\mathfrak A}
\newcommand{ \Bfr}{\mathfrak B}
\newcommand{\Cfr}{\mathfrak C}
\newcommand{\Dfr}{\mathfrak D}

\newcommand{\baP}{\bar P}
\newcommand{\baG}{\bar G}
\newcommand{\baf}{\bar f}
\newcommand{\baB}{\bar B}
\newcommand{\Rep}{\mathcal R}

\title{Transfert du pseudo-coefficient de Kottwitz et formules de caractère
pour la série discrète de ${\rm GL}(N)$ d'un corps local}
\author{P. Broussous \\
                       \\
      Universit\'e de Poitiers \\
       Laboratoire de Mathématiques et Applications\\
       UMR 7348 du CNRS}

\date{\today} 

\begin{document}
\maketitle

\abstract{Let $G$ be the group ${\rm GL}(N,F)$, where $F$ is a non-archimedean locally compact field. 
Using Bushnell and Kutzko's simple types, as well as an original idea of Henniart's, we construct explicit
pseudo-coefficients for the discrete series representations of $G$. As an application we deduce new
formulas for the value of the Harish-Chandra character of certain such representations at certain elliptic regular
elements.}  
\vskip2cm

\centerline{\bf \large  Introduction}
\vskip1cm

  Soit $F$ un corps local non archim\'edien de caract\'eristique
  quelconque et $N>1$ un entier naturel. L'ojectif de cet article est
  d'\'etablir des formules donnant le caract\`ere de certaines
  repr\'esentations de la s\'erie discr\`ete de ${\rm GL}(N,F)$ en certains
  \'el\'ements elliptiques r\'eguliers. Ces formules se veulent \`a la fois
  simples,  et explicites,  au sens o\`u la valeur du caract\`ere est reli\'ee au
  type de Bushnell et Kutzko [BK] de la repr\'esentation.

  Pour arriver \`a ces fins, nous nous servons de la formule de
  Kazhdan reliant la
  valeur du caract\`ere d'une repr\'esentation de la s\'erie discr\`ete \`a
  l'int\'egrale orbitale d'un pseudo-coefficient de cette repr\'esentation
  ([Ka] et [Ba]). L'obtention d'un pseudo-coefficient explicite est
  bas\'ee sur une id\'ee originale de Guy Henniart. Celle-ci consiste \`a
  remarquer les faits suivants.
\medskip

 -- Le pseudo-coefficient construit par Kottwitz [Kott] 
d'une repr\'esentation de carr\'e int\'egrable
 Iwahori-sph\'erique vit dans l'alg\`ebre de Hecke-Iwahori de son type.
\medskip

 -- L'alg\`ebre de Hecke du type d'une repr\'esentation de carr\'e
 int\'egrable de ${\rm GL}(N,F)$ est isomorphe,  via les isomorphismes
 d'alg\`ebres de Hecke de [BK],  \`a une alg\`ebre de Hecke-Iwahori d'un
 autre groupe lin\'eaire $H$ (sur un corps diff\'erent).
\medskip

 L'id\'ee d'Henniart est alors de transf\'erer les pseudo-coefficients des
 repr\'esentations Iwahori-sph\'eriques de la s\'erie discr\`ete de $H$ via
 les divers isomorphismes d'alg\`ebres de Hecke de [BK] en des fonctions sur
 $G$. Les fonctions obtenues sont alors  des candidats pour \^etre des
 pseudo-coefficients des repr\'esentations de la s\'erie discr\`ete de $G$.
Il faut en r\'ealit\'e travailler avec des alg\`ebres de Hecke {\it \`a
  caract\`ere central}, difficult\'e que nous omettons dans cette introduction.  

 \medskip

 Dans ce travail nous montrons que les candidats d'Henniart sont
 effectivement des pseudo-coefficients. Nous calculons leurs
 int\'egrales orbitales pour obtenir des formules de caract\`ere dans deux
 cas :
\medskip

(i)  L'extension de corps $E/F$ qui param\`etrise le type de Bushnell et
Kutzko est non ramifi\'ee et l'\'el\'ement pour lequel on calcule le
caract\`ere est minimal au sens de [BK] (1.4.14),  et engendre dans ${\rm
  M}(N,F)$ une extension non ramifi\'ee.
\smallskip

 (ii) L'extension de corps $E/F$ est totalement ramifi\'ee,  et l'\'el\'ement
elliptique r\'egulier est minimal et engendre une extension totalement
ramifi\'ee. 
\medskip

 Lorsqu'une repr\'esentation cuspidale irr\'eductible $\pi$ d'un groupe
 r\'eductif $p$-adique $G$ est donn\'ee comme une induite compacte
 $\ds \text{c-Ind}{}_K^G\rho$, o\`u $\rho$ est une repr\'esentation
 lisse d'un sous-groupe ouvert compact modulo le centre $K$, le
 caract\`ere de $\pi$ en un \'el\'ement (elliptique) r\'egulier est
 donn\'e par une formule \`a la Frobenius faisant intervenir le
 caract\`ere de $\rho$, comme dans le cas des
 groupes finis (voir par exemple [BH4], Th\'eor\`eme (A.14)). Si $\pi$
 est une repr\'esentation de la s\'erie discr\`ete, elle poss\`ede un
 type $(K,\rho )$, mais ne s'exprime plus comme induite compacte \`a
 partir de ce type si elle n'est pas cuspidale. Il est remarquable
 cependant  que les formules (i)
 et (ii) expriment le caract\`ere de $\pi$ par des formules \`a la
 Frobenius qui utilisent le caract\`ere d'une autre repr\'esentation
 $(K' ,\rho ')$, qui se d\'eduit de $(K,\rho )$ par une modification
 relativement simple. 

\bigskip

 La motivation \`a la base de la r\'edaction de ce travail est
 double. Tout d'abord tr\`es peu de formules sont connues sur le
 caract\`ere d'une repr\'esentation de carr\'e int\'egrable de ${\rm
   GL}(N,F)$ non cuspidale (voir plus bas pour un tr\`es bref rappel
 historique). Ensuite Bushnell et Henniart ont dans une  s\'erie
 d'articles ([BH2] \`a [BH4]) explicit\'e une grande partie de la
 correspondance de Jacquet-Langlands entre s\'erie discr\`ete de ${\rm
   GL}(N)$ et celle d'une de ses formes int\'erieures ${\rm GL}(m,D)$,
 $D$ $F$-alg\`ebre \`a division centrale. Leur explicitation demande
 de conna\^itre le caract\`ere d'une repr\'esentation de la s\'erie
 discr\`ete en suffisament d'\'el\'ements elliptiques
 r\'eguliers. Seules les images par la correspondence de
 certaines repr\'esentations cuspidales sont d\'etermin\'ees. Les
 formules (i) et (ii) devraient permettre d'\'etendre les travaux de
 Bushnell et Henniart \`a certaines repr\'esentations non cuspidales. 
\bigskip

 Le calcul explicite des caract\`eres de repr\'esentations lisses
 irr\'eductibles des groupes r\'eductifs sur $F$ a une longue histoire qui
 remonte au moins aux travaux de Gelfand et Graev [GG]. Pour un tr\`es
 joli historique, nous renvoyons renvoyons \`a l'article de P.J. Sally
 Jr et L. Spice [SaSp]. Cependant les articles de Bushnell et Henniart
 sur l'explicitaton de la correspondance de Jacquet-Langlands ([BH2] \`a [BH4])
 contiennent nombre de formules de caract\`ere in\'edites pour les
 cuspidales, dont il n'est pas fait mention dans [SaSp]. 
\medskip

 Les progr\`es effectu\'es ces cinquante derni\`ere ann\'ees ne concernent
 principalement que les repr\'esentations cuspidales ou les induites
 paraboliques irr\'eductibles. On sait tr\`es peu de choses sur le
 caract\`eres des repr\'esentations des repr\'esentations non cuspidales  de la s\'erie
 discr\`ete (de ${\rm GL}(n,F)$ ou bien d'autres groupes). 
A ma connaissance, seules sont connues
 les valeurs du caract\`ere de la repr\'esentation de Steinberg pour un
 groupe r\'eductif connexe quelconque, et les valeurs du caract\`eres des
 repr\'esentations de niveau $0$ de la s\'erie discr\`ete de ${\rm GL}(m,D)$ en certains
 \'el\'ements minimaux, par les  travaux de Silberger et Zink
 [SZ2]. Notons  que Schneider et Stuhler ont
 obtenu des formules de caract\`eres,  dans une situation
tr\`es g\'en\'erale,  par des consid\'erations homologiques
sur l'immeuble affine de Bruhat-Tits. Cependant leurs formules ne sont
pas exploitables pour obtenir des formules explicites en fonctions des
types, sauf en niveau $0$, o\`u ce travail reste à \'ecrire. 
\bigskip

 Lorsque $E/F$ est non ramifi\'ee, et pour le groupe ${\rm GL}(N,F)$, 
 notre formule g\'en\'eralise celle de
 Silberger et Zink. Nous pensons que le principe du transfert d'un
 pseudo-coefficient devrait permettre d'obtenir des formules dans des
 situations  beaucoup plus g\'en\'erales (en particulier pour d'autres groupes
 r\'eductifs que ${\rm GL}(N)$). 
\bigskip

 Cet article est structuré de  la façon suivante. Le {\S}1 consiste en
 un {\it très bref} rappel de la construction des {\it types simples}
 (de niveau >0), c'est-à-dire des types contenus dans des
 représentations de la série discrète. Un type simple étant fixé, les
 représentations de la série discrète qui le contiennent sont
 construites au {\S}5. 

La structure des algèbres de Hecke-Iwahori est
 rappelée au {\S}2, tandis que les isomorphismes d'algèbres de Hecke
 de Bushnell et Kutzko sont décrits au {\S}3. Le paragraphe {\S}4
 consiste en des définitions et lemmes techniques sur les algèbres de
 Hecke à caractère central. Nous avons eu  besoin d'expliciter
 l'action d'une algèbre de Hecke sphérique sur certains
 modules. N'ayant pas trouvé de références pratiques, nous avons
 préféré donner des démonstrations dans l'annexe A. 

Au {\S}6 nous démontrons que le transfert d'un pseudo-coefficient
satisfaisant à des hypothèses raisonnables est
un pseudo-coefficient. Nous avons besoin pour cela de savoir qu'un
isomorphisme unitaire d'algèbres de Hecke préserve les représentations
{\it tempérées}, fait qui nous a  été communiqué par Bushnell et Henniart
et qui est démontré en Annexe B.  La construction du pseudo-coefficient de
Kottwitz est  rappelée au {\S}7. Dans les {\S}8 et 11 nous
transférons ce pseudo-coefficient dans un premier cas, où $E/F$ est
non ramifié, afin d'obtenir une première formule de caractère. Nous
aurons besoin pour cela de lemmes techniques concernant les
représentations du groupe linéaire sur un corps fini ({\S}9 et
10). Ceux-ci permettent de calculer le caractère d'une
sous-représentation irréductible d'une induite parabolique réductible
en fonction d'un idempotent d'algèbre de Hecke. Malgré leur
simplicité, nous ne leur connaissons pas de références et avons
préféré donner des démonstrations complètes. Enfin dans le {\S}12, le
transfert du pseudo-coefficient de Kottwitz est appliqué dans une
deuxième situation, où $E/F$ est ramifié, pour obtenir une seconde
formule de caractère.

\bigskip

 D'un point de vue technique, ce travail est basée sur les
 isomorphismes d'algèbres de Hecke de [BK]. Nous aurons donc besoin de
 la plupart des concepts et résultats de cette monographie. Pour que
 cet article garde une taille raisonnable, nous supposerons le lecteur
 familier avec le formalisme de Bushnell et Kutzko (strates,
 caractères et types simples). 
 \bigskip

 L'auteur remercie chaleureusement Guy Henniart et Colin Bushnell pour
 de nombreux échanges instructifs et motivants.

\setcounter{section}{-1}

\section{Notations}

 Nous fixons pour tout l'article un corps localement compact, non discret et non
 archimédien $F$. Si $K$ est un tel corps, on note
\medskip

 -- $\ofr_K$ son anneau d'entiers,
\smallskip

 -- $\pfr_K$ l'idéal maximal de $\ofr_K$,
\smallskip

 -- $k_K$ le corps résiduel de $K$,
\smallskip

 -- $\varpi_K$ le choix d'une uniformisante de $K$.
\medskip

 Si $L/K$ est une extension algébrique finie de tels corps, on note
 respectivement $e(L/K)$ et $f(L/K)$ son indice de ramification et son degré
 d'inertie.
\medskip

 On fixe un entier $N\geqslant 2$ ainsi qu'un $F$-espace vectoriel $V$
 de dimension $N$. On note $G={\rm Aut}_F \, (V)\simeq {\rm GL}(N,F)$
 le groupe des automorphismes linéaires de $V$. C'est un groupe
 localement compact et totalement discontinu. Les représentations de
 $G$ considérées dans ce travail seront toujours supposées {\it
   lisses} et {\it complexes}
\medskip

 Si $\mathcal A$ est un anneau unitaire, on note ${\mathcal
   A}^\times$, ou encore $U({\mathcal A})$, son groupe d'unités. En
 particulier, en posant $A={\rm End}_F \, (V)\simeq {\rm M}(N,F)$, on
 a $G=A^\times$. On note ${\rm Rad}\, ({\mathcal A})$ le radical de
 Jacobson de $\mathcal A$.
\medskip

 Soit  $M$  une $K$-algèbre centrale simple sur un corps local
 non-archimédien $K$,  isomorphe à ${\rm M}(n, K)$ pour un  entier $n >0$,
 et soit $\Afr$ un $\ofr_K$-ordre héréditaire de $M$. 
On note $e(\Afr /\ofr_K )$ l'indice de ramification défini par 
$$
\pfr_K \Afr= {\rm Rad}\, (\Afr )^{e(\Afr /\ofr_K )}\ .
$$

 \noi Le groupe des $1$-unités de $\Afr$, défini par $U^1 (\Afr
 )=1+{\rm Rad}\, (\Afr )$ est distingué dans $U(\Afr )$. Si de plus
 $\Afr$ est {\it principal} ({\it cf.} [BK]{\S}(1.1)), le quotient
 $U(\Afr )/U^1 (\Afr )$ est (non canoniquement) isomorphe à 
$$
{\rm   GL}((n/e),k_K )^{\times e(\Afr /\ofr_K )} \ 
$$
\noi On note ${\mathcal K}(\Afr )$ le normalisateur de $\Afr$ dans $M^\times$
: 
$$
{\mathcal K}(\Afr ) =\{ x\in M^\times\ ; \ x\Afr x^{-1} =\Afr \}\ .
$$
\noi Supposons de plus que $E/K \subset M$ est une extension de $K$ et
un sous-anneau unitaire de $M$, de sorte que le centralisateur $B$ de
$K$ dans $M$ est isomorphe à ${\rm M}(n/[E:F], E)$. Alors si $\Bfr$
est un $\ofr_E$-ordre héréditaire de $B$, il existe un unique $\ofr_K$-ordre
héréditaire $\Afr$ de $M$ tel que 
$$
{\mathcal K}(\Bfr )\subset {\mathcal K}(\Afr )\ ,
$$
\noi et de plus $e(\Bfr /\ofr_E ) = e(\Afr /\ofr_K )/e(E/F)$. Pour
plus de détails, nous renvoyons le lecteur à [BK]{\S}(1.2).

\section{Les donn\'ees typiques}

 Puisque le caractère de
Harish-Chandra d'une  représentation irréductible de carré intégrable
 a déjà été étudié par
Silberger et Zink en niveau $0$ [SZ2], nous ne considérerons dans ce
travail que des représentations de niveau $>0$, c'est-à-dire ne
possédant pas de vecteur fixe non nul sous le premier sous-groupe de
congruence
$$
K^1 =1+\varpi_F {\rm M}(N,\ofr_F )\ .
$$

Par le Corollaire (8.5.11), page 304, de [BK], toute représentation
lisse irréductible de carré intégrable de $G$ contient par restriction
un type simple\footnote{Pour une introduction à la théorie générale
  des types, nous renvoyons le lecteur à [BK2]}
 de $G$ au sens de {\it loc. cit.} {\S}5.
\medskip

 Un  type simple $(J,\lambda )$ de niveau $>0$ de $G$ est 
 associ\'e aux  donn\'ees suivantes.

\bigskip

 (a)  Une  strate simple $[\Afr ,n,0,\beta ]$ ([BK](1.5.5), page 43).
 En particulier $\Afr$ est un $\ofr_F$-ordre
 héréditaire principal de $A$, $\beta$ un élément non nul de $A$
 engendrant un corps $E=F[\beta ]$, et normalisant l'ordre $\Afr$. On
 a de plus $\beta\Afr ={\rm Rad}\, (\Afr )^{-n}$.

 On note $B={\rm End}_E (V)$, $\Bfr =\Afr
 \cap B$ ; c'est un $\ofr_E$-ordre héréditaire de $B$ vérifiant
 ${\mathcal K}(\Bfr )\subset {\mathcal K}(\Afr )$. On pose :
$$
e = e(\Bfr /\ofr_E  )=e(\Afr /\ofr_F )/e(E/F) \ .
$$

\smallskip

 (b)  Un caract\'ere simple $\theta\in {\mathcal C} (\Afr ,0,\beta )$
 du groupe $H^1 (\beta ,\Afr )\subset U(\Afr )$ ([BK]{\S}3).
\smallskip

 (c) Une extension $\eta =\eta (\theta )$, dite de Heisenberg,
 de $\theta$ au
 groupe $J^1 (\beta ,\Afr )$ ([BK]{\S}3 et Proposition (5.1.1), page
 158). On a $H^1 (\beta ,\Afr )\subset  J^1 (\beta ,\Afr )\subset
 U(\Afr )$ et $\eta$ est à isomorphisme près l'unique représentation
 lisse iréductible de $J^1 (\beta ,\Afr )$ contenant $\theta$ par
 restriction. 
\smallskip

 (d)  Une $\beta$-extension $\kappa$  de $\eta$ au groupe $J= J(\beta
 ,\Afr )$ ([BK]{\S}3 et Definition (5.2.1), page 166).
\medskip

 Le groupe  $J(\beta ,\Afr )/J^1 (\beta ,\Afr )$ s'identifie canoniquement
 à $U(\Bfr )/U^1 (\Bfr )$, et non canoniquement à 
$$
{\rm GL}(\frac{N}{[E:F]e}, k_E )^{\times e}\ .
$$

 (e) Il existe alors une représentations cuspidale irréductible
$\sigma_0$ de

\noi  ${\rm GL}(N/([E:F]e), k_E )$, telle que 
$$
\lambda= \kappa \otimes \sigma\ ,
$$
\noi où $\sigma= \sigma_0^{\otimes e}$ est vue comme une
représentation de $J(\beta ,\Afr )$, triviale sur $J^1 (\beta ,\Afr
)$. 
\medskip

 Soit ${\mathcal R}_{(J,\lambda )}(G)$ le block  de Bernstein formé
 des représentations lisses de $G$ qui sont engendrées par leur
 composante $\lambda$-isotypique. Alors ([BK] Theorem (7.7.1) et
 Corollary (8.5.11))  ${\mathcal R}_{(J,\lambda )}(G)$ contient toujours
 des représentations irréductibles de carré intégrable.

\bigskip

 Nous aurons besoin aussi de travailler avec un  type modifi\'e 
$(J' , \lambda ')$ défini dans [BK]{\S}(5.6), page 189, type pour la
 même composante de $G$ que $(J,\lambda )$. Il est construit de la
 façon suivante.
\bigskip

\smallskip

 On fixe une extension non ramifi\'ee $K/E$ de degr\'e $f:=N/([E:F]e)$ telle que
 $K^\times \subset {\mathcal K}(\Bfr )$. L'ordre $\Cfr =\Bfr \cap C =
 \Afr \cap C$ est un ordre héréditaire {\it minimal} (ou d'Iwahori) de
 $C$ et ${\mathcal K}(\Cfr )\subset {\mathcal K}(\Bfr )\subset {\mathcal
   K}(\Afr )$. 
\smallskip

 On fixe un ordre maximal $\Cfr_M \supset \Cfr$ et on note $\Afr_M$
 et $\Bfr_M$ les ordres correspondants dans $A$ et $B$
 respectivement. En d'autres termes, $\Afr_M$ (resp. $\Bfr_M $) est
 l'unique ordre héréditaire de l'algèbre $A$ (resp. de l'algèbre $B$)
 tel que ${\mathcal K}(\Cfr_M )\subset {\mathcal K}(\Afr_M )$ (resp. tel
 que ${\mathcal K}(\Cfr_M )\subset {\mathcal K}(\Bfr_M )$). 
L'ordre $\Bfr_M$ est maximal, tandis que l'ordre $\Afr_M$
 ne l'est pas toujours.
\smallskip

  On note $\theta_M$ le transfert de $\theta$ \`a $H^1 (\beta
 ,\Afr_M )$ (voir [BK]{\S}(3.6) pour la définition du transfert d'un
  caractère simple entre deux ordres). Soit   $\eta_M$ l'extension de 
Heisenberg de $\theta_M$, et $\kappa_M$ la
 $\beta$-extension de $\eta_M$ attach\'ee \`a $\kappa$ comme dans
 [BK],  page 167, Th\'eor\`eme (5.2.3). C'est une représentation de $J
 (\beta ,\Afr_M )$.  On  restreint $\kappa_M$ en
  une repr\'esentation de
 $J' = U(\Bfr )J_M^1$, o\`u $J_M^1 =J^1 (\beta , \Afr_M )$.
\smallskip

  Alors $\lambda '$ est la repr\'esentation de $J'$ donn\'ee par
 $\kappa_M\otimes \sigma$. On voit ici  $\sigma$ comme une
 repr\'esentation de $J'$ via l'isomorphisme canonique
$$
J' /U^1 (\Bfr ) J^1_M \simeq  U(\Bfr )/U^1 (\Bfr )\ .
$$

\bigskip

 Le fait que la  paire $(J' , \lambda ')$ est un type pour la 
m\^eme composante de  Bernstein que $(J,\lambda )$ vient de ce que ces
 deux repr\'esentations
 s'induisent en des repr\'esentations irr\'eductibles \'equivalentes
 de $U(\Bfr )U^1 (\Afr )$ ([BK], Prop. (5.5.13), page 185).

\section{L'alg\`ebre de Hecke-Iwahori}

 Posons $H={\rm GL}_K (V)= C^{\times}$. Alors $U(\Cfr )$ est un
 sous-groupe d'Iwahori de $H$, que l'on notera parfois $I$. 
On fixe une base $(v_1 ,...,v_e  )$ de
 $V$ sur $K$, de sorte que dans l'identification $H\simeq {\rm
   GL}(e,K)$, $U(\Cfr)$ soit form\'e des matrices de $GL (d,\ofr_K )$
 qui sont triangulaires sup\'erieures modulo $\pfr_K$.  On supposera
alors que  $\Cfr_M$ est l'ordre maximal standard ${\rm M}(e,\ofr_K )$,
de sorte que $U(\Cfr_M )$ est le compact maximal standard de $H$.

\medskip

 On note $\HH_0$ l'alg\`ebre de Hecke-Iwahori $\HH (H//U(\Cfr ))\simeq
 \HH (H,{\mathbf 1}_{U(\Cfr )})$.  Le produit de convolution est ici
 d\'efini en choisissant comme mesure de Haar sur $H$ {\it celle qui
   donne un  volume $1$ \`a $U(\Cfr )$}. Elle admet pour base les
 fonctions caract\'eristiques

$$
 f^0_w = {\mathbf 1}_{U(\Cfr )w U(\Cfr )}
$$
\noi o\`u $w$ d\'ecrit $W_0^{\rm aff }$ le groupe de Weyl affine
\'etendu de $H$. On note $W_0 = W_0^{\rm sph}$ le groupe de Weyl
sph\'erique de $H$ relativement au sous-groupe de Borel des matrices
triangulaires sup\'erieures de $H$, et $S\subset W_0$ le syst\`eme
d'involutions habituel engendrant $W_0$~:

$$
S =\{ s_i \ ; \ i=1,...,e-1\}
$$
\noi o\`u $s_i$ est la matrice correspondant \`a la transposition
$(i,i+1)$. Rappelons que $(W_0 ,S)$ est un groupe de Coxeter. Le
sous-groupe $U(\Cfr )$ de $U(\Cfr_M )$ est le groupe "$B$" d'une
 ${\rm BN}$-paire (ou syst\`eme de Tits)
du groupe $U(\Cfr_M )$,  de groupe de Coxeter $(W_0 ,S )$.

 La sous-alg\`ebre $\HH_0^{\rm sph} =  \HH (U(\Cfr_M )//U(\Cfr ))
 \subset \HH_0$ est  form\'e des fonctions \`a support dans $U(\Cfr_M
 )$  Elle admet pour base les $f_w^0$, $w\in W_0$.
\bigskip

 Nous aurons besoin pour la suite d'un syst\`eme de repr\'esentants
 $\mathcal P$ des $H$-classes de conjugaison des sous-groupes
 parahoriques de $H$. Un sous-groupe parahorique de $H$ est de la
 forme $U(\Dfr)$, o\`u $\Dfr$ est un $\ofr_K$-ordre h\'er\'editaire de
 ${\rm End}_K (V)$. Dans la classe de conjugaison d'un parahorique,
il y a toujours un $U(\Dfr )$ tel que $\Cfr \subset \Dfr \subset
\Cfr_M$, de sorte que par propri\'et\'e des syst\`emes de Tits, on a

$$
U(\Dfr )=U(\Cfr )\langle T\rangle U(\Cfr )
$$
 \noi pour un certain $T\subset S$.  On note que deux parahoriques
$$
U(\Dfr_i ) =U(\Cfr )\langle T_i \rangle U(\Cfr )\  , \ T_i \subset S\ , \ i=1, 2
$$
\noi sont conjugu\'es si, et seulement si, $T_1$ et $T_2$ sont
conjugu\'es sous l'action du groupe $\langle \Pi \rangle$, o\`u

$$
\Pi =
\left(
\begin{array}{cccccc}
0              & 1  & 0  & 0  & \cdots &  0 \\
0              & 0  & 1  & 0  & \cdots  & 0 \\
\vdots             &     &  &  &  & \vdots\\
0              & \cdots & \cdots & 0 & 1& 0 \\
0              & \cdots & \cdots & \cdots & 0 & 1 \\
\varpi_K  & 0  & \cdots & \cdots  & \cdots & 0
\end{array}
\right)
$$

 Nous fixons

$$
{\mathcal P}=\{ P_T= U(\Cfr )\langle T\rangle   U(\Cfr ) \   , \ T\in \Theta\  \}
$$
\noi o\`u $\Theta$ d\'esigne un syst\`eme de
repr\'esentants des orbites de $\langle \Pi \rangle$ dans l'ensemble
des parties de $S$ (ce qui est
un abus de langage car l'action de $\Pi$ par conjugaison ne laisse pas
$S$ stable).

\section{Les isomorphismes d'alg\`ebres de Hecke}

Nous nous r\'ef\'erons ici \`a [BK] (5.6), eux-m\^emes se r\'ef\'erant
 \`a  [HM]. 
Le but de cette section est de rappeler la forme
 explicite de l'isomorphisme d'alg\`ebres de Hecke :
$$
\Psi~: \HH_0 =\HH (H,{\mathbf 1}_{U(\Cfr )})\lra 
 \HH_{\lambda '}=\HH (G,\lambda ')\  ,
$$
\noi o\`u du moins, c'est ce qui nous int\'eresse dans un premier 
temps, de sa restriction :
$$
\Psi \ : \HH_0^{\rm sph} =\HH (U(\Cfr_M )//U(\Cfr ))\lra \HH_{\lambda '}\
$$
\noi Nous allons rappeler sa construction.
\medskip

 Rappelons que $J_M = U(\Bfr_M )J_M^1$, $J' =U(\Bfr )J_M^1$, et que le type $(J',\lambda ')$ est
 donn\'e par $(\kappa_M )_{\vert J'}\otimes \sigma$, o\`u :
\medskip

-- $(J_M, \kappa_M )$ est la $\beta$-extension attach\'ee \`a $(J, \kappa)$,
\smallskip

 --  $\sigma$ est vue comme une repr\'esentation de $U(\Bfr )J_M^1$ 
triviale sur $U^1 (\Bfr )J_M^1$,
 via l'isomorphisme canonique
$$
U(\Bfr )J_M^1 /U^1 (\Bfr )J_M^1 \simeq U(\Bfr )/U^1 (\Bfr )\ .
$$

 Le quotient $J_M /J_M^1 =: {\bar G}$ est canoniquement isomorphe \`a 
$U(\Bfr_M )/U^1 (\Bfr_M ) \simeq {\rm GL}(ef,k_E )$.
  Le sous-groupe $J'/J_M^1$ s'identifie alors \`a un sous-groupe
 parabolique $\bar P$ de $\bar G$,
 de sous-groupe de Levi ${\bar L} =U(\Bfr )/U^1 (\Bfr )\simeq {\rm GL}(f,k_E )^e$, et
 de radical unipotent $ {\bar U}= U^1 (\Bfr )/U^1 (\Bfr_M )$.

 On peut donc regarder $\sigma$ comme une repr\'esentation de $\bar L$ et consid\'erer
 son inflation \`a $\bar P$, que l'on note par le m\^eme symbole.

 Rappelons le r\'esultat principal de [HM].

\begin{proposition}  Il existe un unique isomorphisme d'alg\`ebres :
$$
\Upsilon~: \HH_0^{\rm sph} \lra \HH ({\bar G}, \sigma_{\vert \bar P})
$$
\noi  qui pr\'eserve (dans un sens que nous allons pr\'eciser plus loin) le support des fonctions.
\end{proposition}

 \begin{remark} La mesure de comptage sur $\bar G$ qui d\'efinit le produit de convolution
 de $ \HH ({\bar G}, \sigma_{\vert \bar P})$ donne le volume $1$ \`a tout singleton de $\bar G$.
\end{remark}

 Notons encore $W_0$ le groupe des matrices de permutation de ${\rm GL}(e,k_E )$ plong\'e par blocs 
dans ${\rm GL}(ef,k_E )$  (un "$0$" est remplac\'e par un bloc $0_f$ et un "$1$" par un bloc $I_f$).
\medskip

 L'entrelacement de $\sigma_{\vert \bar P}$ dans $\bar G$ est alors donn\'e par ${\bar P}W_0 {\bar P}$.
 Notons $X_0$ l'espace de $\sigma$ et $X=X_0^{\otimes e}$ celui de $\sigma$. Pour chaque $w\in W_0$,
 on note ${\bar f}_w$, l'\'el\'ement de $\HH ({\bar G}, \sigma_{\vert \bar P})$, de support ${\bar P}W_0 {\bar P}$, donn\'e par
$$
{\bar f}_w (p_1 wp_2 )=
\frac{1}{\vert {\bar P}\vert}
 {\tilde \sigma}(p_1 )\circ T_w \circ {\tilde \sigma }(p_2 )\ ,
\ p_1 , \ p_2 \in {\bar P}\ ,
$$
\noi o\`u $T_w$ d\'esigne l'endomorphisme de $\tilde X$ donn\'e par
$$
T_w (x_1 \otimes x_2 \otimes \cdots \otimes x_e )=x_{w(1)}\otimes x_{w(2)}\otimes \cdots x_{w(e)}\ .
$$

Ces ${\bar f}_w$, $w\in W_0$, forment une base de
 $\HH ({\bar G}, \sigma_{\vert \bar P})$.

\bigskip

Notons $\bar B$ le sous-groupe de Borel standard sup\'erieur du groupe ${\rm GL}(e,k_K )$.  Pour $w\in W_0$, 
 notons  $\baf_w^0$ l'\'el\'ement de l'alg\`ebre de Hecke sph\'erique $\HH ({\rm GL}(e,k_K ),{\bar B} )$ donn\'e par
$$
\baf_w^0 =\frac{1}{\vert {\bar B}\vert}\ {\mathbf 1}_{{\bar B}w{\bar B}}\ .
$$
\noi Les $\baf_w^0$, $w\in W_0$, forment une base de cette alg\`ebre. On a alors un isomorphisme d'alg\`ebres ([HM]) :
$$
\HH (\baG ,\sigma_{\vert \baP})\lra \HH ({\rm GL}(e,k_K ), {\bar B})
$$
\noi donn\'e sur la base par :
$$
\baf_w \mapsto \baf_w^0\ .
$$

Quant \`a elle, l'alg\`ebre $\HH ({\rm GL}(e, k_K ),\baB )$ est
naturellement isomorphe \`a $\HH_0^{\rm sph} (H//U(\Cfr )) =\HH
(U(\Cfr_M )//U(\Cfr ))$ par
$$
\baf_w^0 \mapsto f_{w}^0 ={\mathbf 1}_{U(\Cfr )wU(\Cfr )} \ , \ w\in
W_0\ .
$$

Finalement notre isomorphisme $\Upsilon$ s'obtient en composant ces
isomorphismes et est donn\'e sur les bases par :

$$\begin{array}{ccccc }
\Upsilon~: & \HH_0^{\rm sph} & \lra & \HH (\baG ,\sigma_{\vert \baP })
&  \\
           &  f_w^0        &\mapsto & \baf_w\ &       \ , \  w\in W_0\\
\end{array}
$$

La derni\`ere \'etape est de plonger $\HH (\baG ,\sigma_{\vert \baP})$
dans l'alg\`ebre $\HH (G,\lambda ')$ via l'isomorphisme d'alg\`ebres
$$
\HH (\baG ,\sigma_{\vert \baP})\simeq \HH (J_M ,\lambda ')\subset \HH
(G,\lambda ')
$$
\noi donn\'e par [BK], Lemma (5.6.3), page 189. D\'ecrivons cet
isomorphisme. Pour cela notons $Y$ l'espace de $\kappa_M$.
\medskip

\noi {\bf N.B.} Je n'ai pas les m\^emes notations que [BK] : les r\^oles
de $X$ et $Y$ sont interchang\'es.
\medskip

\`A $\varphi \in \HH (\baG ,\sigma_{\vert \baP})$, il associe la
fonction $\varphi '$~: $J_M\lra {\rm End}_{\CC}({\tilde X}\otimes
{\tilde Y})$ donn\'ee par
$$
\varphi '(g)={\tilde \kappa}_M  (g)\otimes \varphi ({\bar g})
$$
\noi o\`u $\bar g$ d\'esigne l'image de $g$ dans $J_M /J_M^1 \simeq {\bar
  G}$.
\medskip

 Notons que, m\^eme si c'est tacite dans [BK],  pour que $\varphi\mapsto
 \varphi '$ soit un morphisme d'alg\`ebres, il nous faut imposer :

\begin{notation} La mesure de Haar $\mu_G$ sur $G$ est fix\'ee de telle
  sorte que $\mu_G (J_M^1 )=1$.
\end{notation}

\medskip

\noi {\bf N.B.} Il y a une petite faute de frappe dans [BK], o\`u il est
\'ecrit $g\in G$, au lieu de $g\in J_M$.
\medskip

 Nous noterons
$$
\Psi_0~: \HH_0^{\rm sph}\lra \HH (J_M ,\lambda ' )\subset \HH (G,\lambda
' )
$$

\noi la composition des divers isomorphismes d'alg\`ebres de Hecke
construits pr\'ec\'edemment.
\bigskip

\begin{theorem} ([BK] Proposition (5.5.11), page 185.) L'entrelacement de $(J' , \lambda ')$ dans $G$ est $J'
  W_{0}^{\rm aff} J'$, o\^u, rappelons le, $W_{0}^{\rm aff}$ est le
  groupe de weyl affine de $H$.
\end{theorem} 

 Consid\'erons l'\'el\'ement de $\HH (H,{\mathbf 1}_{U(\Cfr )})$ donn\'e par
 $\zeta ={\mathbf 1}_{U(\Cfr )\Pi U (\Cfr )}={\mathbf 1}_{\Pi U(\Cfr
   )}$.

\begin{theorem} ([BK] Main Theorem (5.6.6), page 190.) (i) Il existe
  un \'el\'ement non nul  $\psi$ de $\HH (G, \lambda ' )$ de support $J'
  \Pi J'$ et celui-ci est unique \`a un scalaire pr\`es.
\smallskip

\noi (ii) Pour tout \'el\'ement $\psi$ comme en (i), il existe un unique
morphisme d'alg\`ebres 
$$
\Psi~: \ \HH (H, {\mathbf 1}_{U(\Cfr )})\lra \HH (G,\lambda ')
$$
\noi qui prolonge $\Psi_0$ et v\'erifie $\Psi (\zeta )=\psi$.
\smallskip

De plus, une telle application $\Psi$ est un isomorphisme d'alg\`ebres
qui pr\'eserve les supports au sens o\`u le support de $\Psi ({\mathbf
  1}_{U(\Cfr )wU(\Cfr )})$ est $J' wJ'$, pour tout $w\in W_{0}^{\rm
  aff}$.
\end{theorem}

 Pour toute la suite de ce travail, nous fixons un tel isomorphisme
 $\Psi$. Nous supposerons :

\begin{hypothese} L'isomorphisme $\Psi$ est unitaire au sens 
de la   d\'efinition (5.6.16) de [BK].
\end{hypothese}

 Rappelons ([BK](4.3), page 153) que l'alg\`ebre $\HH (G,\lambda ')$ est
 canoniquement munie d'une (anti-)involution (semi-lin\'eaire) canonique 
 $\Phi\mapsto {\bar \Phi}$. L'isomorphisme $\Psi$ est alors dit {\it
   unitaire} si, avec les notations pr\'ec\'edentes, on a $\psi \star
 {\bar \psi}=1$. Notons que cette condition d\'etermine $\Psi$  un
 facteur pr\`es dans $\{ u\in \CC\ ; \ \vert u\vert =1\}$.

\section{Alg\`ebres de Hecke \`a caract\`ere central}

 Notons $\omega_0$ un caract\`ere  de $K^\times$ trivial sur
 $\ofr_K^\times$   et $\omega_0
 {\mathbf 1}_{U(\Cfr )}$ le caract\`ere de $K^{\times}U(\Cfr
 )$ prolongeant $\omega_0$ et ${\mathbf 1}_{U(\Cfr )}$. 
Notons aussi ${\mathcal H}(H,\omega_0 {\mathbf 1}_{U(\Cfr )})$
 l'alg\`ebre de Hecke des fonctions $f$~: $H\lra \CC$ qui
\medskip

 -- sont bi-invariantes sous l'action de l'Iwahori ${U(\Cfr )}$,
\smallskip

 -- v\'erifient $f(zg)=\omega_0^{-1}(z)f(g)$, $g\in H$, $z\in
 K^\times$,
\smallskip

 -- sont \`a support compact modulo le centre $K^\times$.
\medskip

 On note $\mu_H$ la mesure de Haar sur $H$ donnant un volume $1$ \`a
 l'Iwahori $U(\Cfr )$ et $\mu_{K^\times}$ la mesure de Haar sur
 $K^\times$ donnant un volume $1$ \`a $\ofr_K^\times$. On d\'esigne alors
 par $\mu_{H/K^\times}$ la mesure quotient d\'efinie par
$$
\int_H f(x) d\mu_H (x)=\int_{H/K^\times} \big\{ \int_{K^\times}
f(zy)d\mu_{K^\times}(z)\big\} \ d\mu_{H/K^\times} (\dot{y} )\ , \ f\in
{\mathcal C}_c^{\infty}(H)\ .
$$

 La structure d'alg\`ebre sur $\HH (H, \omega_0 {\mathbf 1}_{U(\Cfr )})$
 est donn\'ee par la convolution :
$$
f_1 \star f_2 (g)=\int_{H/K^\times}f_1 (y)f_2
(y^{-1}g)d\mu_{H/K^\times}(\dot{y})\ , \ f_1 ,\ f_2 \in  \HH (H,
\omega_0 {\mathbf 1}_{U(\Cfr )})\ .
$$

 Soit $(\sigma ,{\mathcal W})$ un objet de ${\mathcal R}_{\rm
   Iw}(H)$,  de caract\`ere central $\omega_0$. Alors le $\CC$-espace
 vectoriel $N=\WW^{U(\Cfr )} =\WW^{\omega_0 {\mathbf 1}_{U(\Cfr )}}$
   est  muni d'une
 structure de $\HH (H,{\mathbf 1}_{U(\Cfr )})$-module \`a gauche via
$$
\sigma (f).w = \int_{H}f(h)\sigma (h).w\, d\mu_H (h)\ , \ f\in\HH
(H,{\mathbf 1}_{U(\Cfr )}) , \ w\in N\ ,
$$
\noi et d'une structure de  $\HH (H, \omega_0 {\mathbf 1}_{U(\Cfr
  )})$-module \`a gauche via~:
$$
\sigma (f).w =\int_{H/K^\times} f(h)\sigma
(h).w\ d\mu_{H/K^\times}(\dot{h})\ ,
$$
\noi $f\in  \HH (H, \omega_0  {\mathbf 1}_{U(\Cfr )})$, $w\in
     {\mathcal W^{\omega_0 {\mathbf 1}_{U(\Cfr )}}}$.

On a un morphisme surjectif d'alg\`ebres :
$$
P_{\omega_0}~: \ \HH (H, {\mathbf 1}_{U(\Cfr )})\lra \HH (H,\omega_0
{\mathbf 1}_{U(\Cfr )})
$$
\noi donn\'e par
$$
P_{\omega_0}(f)(y)=\int_{K^{\times}}\omega_0 (z)
f(zy)d\mu_{K^\times}(z) =\sum_{n=-\infty}^{+\infty}\omega_0 (\varpi_K
)^n f(\varpi_K^n y)\ ,
$$
\noi $y\in H$,  $f\in \HH (H,{\mathbf 1}_{U(\Cfr )})$ (voir
     [BK](6.1.6), page 201). 

\begin{lemma} Soit $(\sigma , {\mathcal W})$ une repr\'esentation de
  ${\mathcal R}_{\rm Iw}$ de caract\`ere central $\omega_0$. 
\medskip

\noi a) Le diagramme suivant est commutatif :
$$
\begin{array}{cccc}
    &   \HH (H,{\mathbf 1}_{U(\Cfr )}) & \rightarrow & {\rm
    End}_{\CC}\,  (N)\\
P_{\omega_0} &  \downarrow & \nearrow &    \\
     & \HH (H,\omega_0 {\mathbf 1}_{U(\Cfr )}) & & 
\end{array}
$$
\noi o\^u les fl\`eches non nomm\'ees d\'ecoulent des structures de module de
$N$. 
\medskip

b) Si $\sigma$ est admissible, alors pour tout $f\in \HH (H,{\mathbf
  1}_{U(\Cfr )})$ on a 
$$
{\rm Tr}(\sigma (f), N)={\rm Tr}(\sigma (P_{\omega_0}(f)),
N)\ .
$$
\end{lemma}

\noi {\it Preuve}. Le point a) d\'ecoule d'un calcul imm\'ediat et b) en
est une cons\'equence \'evidente.
\bigskip

De fa\c con similaire, nous allons d\'efinir une alg\`ebre $\HH
(G,\lambda'\omega )$ comme dans [BK](6.1), page 199. Puisque la
repr\'esentation $(J' ,\lambda ')$ est irr\'eductible, le lemme de Schur
affirme qu'il existe un caract\`ere $\omega_\lambda $ de
$\ofr_F^{\times}$ tel que $\lambda'_{\mid \ofr_F^{\times}}$ est un multiple
  de $\omega_\lambda$. Fixons un caract\`ere $\omega$ de $F^\times$ tel
  que $\omega_{\mid \ofr_F^\times} =\omega_\lambda$. 

 On d\'efinit l'alg\`ebre de Hecke $\HH (G,\lambda '\omega )$ comme
 l'ensemble des fonctions $f$~: $G\lra {\rm End}_\CC
 (\check{W}_\lambda' )$ v\'erifiant~:
\medskip

 -- $f(j_1 gj_2 )= (\lambda '\omega )^{\vee}(j_1 )f(g) (\lambda
 '\omega )^{\vee} (j_2 )$, $j_1$, $j_2\in F^\times J'$, $g\in G$
\smallskip

 -- $f$ est \`a support compact modulo le centre  $F^{\times}$.
\medskip

 On fixe une mesure de Haar $\mu_{F^{\times}}$ sur $F^{\times}$ de
 sorte que $\mu_{F^{\times}}(\ofr_F^\times )=1$ et on  note
 $\mu_{G/F^{\times}}$ la mesure quotient correspondante sur
   $G/F^{\times}$. C'est cette derni\`ere mesure qui nous permet de
   d\'efinir un produit de convolution sur $\HH (G,\lambda '\omega )$. 
\medskip

 On a aussi un morphisme surjectif d'alg\`ebres 
$$
{P}_{\omega}~:\ \HH (G,\lambda ')\lra \HH (G, \lambda '\omega )
$$
\noi donn\'e par
$$
P_{\omega}(f)(g)=\int_{F^\times}\omega (z)f(zg)\ d\mu_{F^{\times}}(z)
= \sum_{n=-\infty}^{+\infty} \omega (\varpi_F )^n f(\varpi_F^n g)
$$
\noi $g\in G$, $f\in \HH (G,\lambda ')$. 
\bigskip

 Consid\'erons une repr\'esentation lisse $(\pi ,\VV )$ de $G$, de
 caract\`ere central $\omega$. Elle donne lieu \`a deux $\CC$-espaces
 vectoriels canoniquement isomorphes :
$$
M={\rm End}_{J'}\, (\lambda ' ,\pi )\text{ et } M_\omega ={\rm
  End}_{F^\times J'}(\lambda '\omega ,\pi )\ .
$$
\noi Plus exactement l'inclusion naturelle $M_{\omega}\subset M$ est
en fait une \'egalit\'e. Comme il est rappel\'e dans l'annexe A, le
$\CC$-espace $M=M_\omega$ est muni d'une structure de $\HH (G,\check{\lambda
'})$-module \`a droite, ainsi que d'une structure de $\HH (G,\check{\lambda
'\omega} )$-module \`a droite. 

On a un morphisme surjectif d'alg\`ebres 
$$
{\tilde P}_\omega ~: \ \HH (G,\check{\lambda '})\lra \HH (G,
\check{\lambda '\omega})
$$
\noi d\'efini par la m\^eme formule que pour $P_\omega$. 

\begin{lemma} 1) On a un diagramme commutatif :
$$
\begin{array}{rccl}
     & \HH (G,\check{\lambda '}) & \lra & {\rm End}_\CC\, (M)\simeq {\rm
    End}_\CC\, (M_\omega )\\
P_{\omega} & \downarrow & \nearrow & \\
  & \HH (G,\check{\lambda '\omega} )
\end{array}
$$
\noi 2) Pour tout $f\in \HH (G,\check{\lambda '} )$, on a 
$$
{\rm Tr}\, (f,M)={\rm Tr}\, (P_\omega (f), M_{\omega})\ .
$$
\end{lemma}

\noi {\it D\'emonstration}. Le point 2) d\'ecoule \'evidemment de 1). 

Soient $f\in \HH (G, \check{\lambda '})$  et $\varphi\in M\simeq
M_\omega$. Notons $f_\omega =P_\omega (f)$. Il s'agit de d\'emontrer que
$\varphi . f = \varphi_\omega .f_\omega$. D'apr\`es les lemmes (A2) et
(A4) de l'annexe A, on a les formules :
$$
\varphi .f=\int_G \pi (x)\circ \varphi \circ f(x^{-1})\, dx
$$
\noi et

\begin{eqnarray}
\varphi . f_\omega &=& \int_{G/F^\times}\pi (x)\circ \varphi \circ f_{\omega}
  (x^{-1})\, d\mu_{G/F^\times} (\dot{x})\\
 &=& \int_{G/F^\times}\big(\,  \pi (x)\circ \varphi \circ \int_{F^\times}\omega
  (z) f(z^{-1}x^{-1})\, d\mu_{F^\times} (z)\, \big) \,
  d\mu_{G/F^\times}(\dot{x})\\
 &=& \int_{G/F^\times} \, \int_{F^\times} \big(\,  \pi (zx)\circ
  \varphi \circ f((zx)^{-1})\, \big) \, d\mu_{F^\times}(z) \,
  d\mu_{G/F^\times}(\dot{x}) \\
 &=& \int_{F}\pi (x)\circ \varphi \circ f(x^{-1})\, dx\\
 &=& \varphi . f
\end{eqnarray}
\noi ce qu'il fallait d\'emontrer.

\bigskip

 Rappelons que les alg\`ebres $\HH (G,\check{\lambda '})$ et $\HH
 (G,\lambda ' )$ sont anti-isomorphes via l'op\'eration $f\mapsto f^*$,
 d\'efinie par
$$
f^* (g) =f(g^{-1})\check{}\ , \ g\in G , \ f\in \HH (G,\check{\lambda
  '})
$$
\noi o\^u on a not\'e $a\check{}$ le transpos\'e dans ${\rm End}_\CC \,
(\check{W})$ d'un \'el\'ement de ${\rm End}_\CC \, (W)$. On d\'efinit de
fa\c con similaire un anti-isomorphisme d'alg\`ebres entre  $\HH
(G,\check{\lambda '\omega})$ et $\HH (G,\lambda '\omega  )$ que l'on
note encore $f\mapsto f^*$. Ceci permet de munir $M$
(resp. $M_\omega$) d'une structure de $\HH (G,\lambda ')$-module \`a
gauche (resp. $\HH (G,\lambda '\omega )$-module \`a gauche) par la
formule~:
$$
f.\varphi = \varphi .f^*\ .
$$

\begin{lemma} 1) On a un diagramme commutatif :
$$
\begin{array}{rccl}
     & \HH (G,{\lambda '}) & \lra & {\rm End}_\CC\, (M)\simeq {\rm
    End}_\CC\, (M_\omega )\\
P_{\omega} & \downarrow & \nearrow & \\
  & \HH (G,{\lambda '\omega} )
\end{array}
$$
\noi 2) Pour tout $f\in \HH (G,{\lambda '} )$, on a 
$$
{\rm Tr}\, (f,M)={\rm Tr}\, (P_\omega (f), M_{\omega})\ .
$$
\end{lemma}

\noi {\it D\'emonstration}. A nouveau 2) d\'ecoule directement de 1). Au
vu du lemme pr\'ec\'edent, pour montrer 1), il suffit d'\'etablir la
commutativit\'e du diagramme suivant :
$$
\begin{array}{lcccr}
&    &  \check{} & & \\
& \HH (G,\check{\lambda '}) & \lra & \HH (G,\lambda ') & \\
{\check P}_\omega & \downarrow &    & \downarrow & P_\omega \\
& \HH (G,\check{\lambda '\omega}) & \lra & \HH (G,\lambda '\omega ) &
\\
 &   &  \check{} &  & 
\end{array}
$$
\noi En effet, si $f\in \HH (G,\check{\lambda '})$ et $g\in G$, on a 
$$
{\tilde P}_\omega (f)(g)=\int_Z \omega (z) f(zg)\, d\mu_Z (z)
$$
\noi et 

\begin{eqnarray}
\big( {\tilde P}_\omega (f)(g^{-1})\big)\check{} &=& \int_{Z}\omega (z)
f(zg^{-1})\check{}\, d\mu_Z (z)\\
       &=& \int_Z \big( \omega (z^{-1})f((z^{-1}g)^{-1})\big)\check \,
d\mu_Z (z)\\
       &=& \int_Z \omega (t) f^{*} (tg)\, d\mu_Z (t)\\
       &=& P_\omega (f^{*})
\end{eqnarray}

\noi ce qu'il fallait d\'emontrer.

\bigskip

Soit $\HH (G,\omega )$ l'espace des fonctions $f$ localement
constantes sur $G$, \`a support compact modulo $F^\times$ et v\'erifiant
$f(zg)=\omega^{-1}(z)f(g)$, $f\in F^\times$, $g\in G$. On le munit
d'une structure d'alg\`ebre via le produit de convolution 
$$
f_1\star f_2 (g)=\int_{G/F^\times} f_1 (x) f_2
(x^{-1}g)\ d\mu_{G/F^\times}(\dot{x})\ .
$$
\noi On note $e_{\lambda '\omega}$ l'idempotent de $\HH (G,\omega )$
d\'efini par 
$$
e_{\lambda '\omega} (x)=
\left\{
\begin{array}{ll}
\mu_{G/F^\times} (J' F^\times )\ {\rm dim}\ (\lambda '){\rm
  Tr}\ (\lambda ' \omega (x^{-1})) & \text{si } x\in J' F^\times \\
0 & \text{sinon}
\end{array}\right.
$$
\noi de sorte que $e_{\lambda ' \omega}\star {\mathcal V}$ est la
composante $\lambda '$-isotypique ${\mathcal V}^{\lambda '}$ 
 de $\mathcal V$, pour toute
repr\'esentation lisse $\mathcal V$ de caract\`ere central $\omega$. 
\medskip

 Par [BK], Proposition (4.2.4) et Remark (4.2.6), on a un isomorphisme
 canonique d'alg\`ebres de Hecke~:
$$
\Upsilon_{\lambda '\omega} ~: \HH (G ,\lambda '\omega )\otimes_\CC {\rm End}_\CC
\ (W_{\lambda '})\lra e_{\lambda '\omega}\star \HH (G,\omega )\star
e_{\lambda '\omega}\ .
$$

Rappelons bri\`evement comment cet isomorphisme est construit. On
commence par identifier ${\rm End}_\CC\ (W_{\lambda '})$ \`a $W_{\lambda
  '}\otimes \check{W}_{\lambda '}$ de la fa\c con canonique
habituelle. Si $\varphi\in \HH (G,\lambda '\omega )$, $w\in W_{\lambda
  '}$, $\check{w}\in \check{W}_{\lambda '}$, l'image $\Phi
=\Upsilon_{\lambda '\omega} ( \varphi \otimes w\otimes\check{w})$ est donn\'ee par 
$$
\Phi (g) = {\rm dim}\ (\lambda ')\  \langle w, \varphi
(g)\check{w}\rangle\ .
$$
En particulier si $(w_i )$ est une base de $W_{\lambda '}$, de base
duale $(\check{w}_i )$, on a $\ds {\rm id}_{W_{\lambda '}} =\sum_i w_i
\otimes \check{w}_i$, et 
$$
\Upsilon_{\lambda '\omega} (\varphi \otimes {\rm id}_{W_{\lambda '}})(g)={\rm
  dim}\ (\lambda ')\sum_i \langle w_i
,\varphi (g)\check{w}_i\rangle = {\rm dim}\ (\lambda '){\rm
  Tr}_{W_{\lambda '}}\ (\varphi (g)) ,
$$

\noi $g\in G\ , \ \varphi \in\HH (G,\lambda '\omega )$. Pour r\'esumer :

\begin{lemma} Pour tout $f\in \HH (G,\lambda '\omega )$, on a 
$$
\Upsilon_{\lambda '\omega} (f\otimes {\rm id}_{W_{\lambda '}})={\rm
  Tr}_{W_{\lambda '}}\circ f\ .
$$
\end{lemma}
   
 Oubliant le caract\`ere central, on peut, comme dans [BK] (4.2),
 consid\'erer l'idempotent  $e_{\lambda '}$ de $\HH (G)$ 
associ\'ee \`a la repr\'esentation $(J',\lambda ' )$ et l'alg\`ebre de Hecke
$e_{\lambda '}\star \HH (G)\star e_{\lambda '}\subset \HH (G)$. On a
alors un isomorphisme canonique
$$
\Upsilon_{\lambda '}~ : \HH (G,\lambda ')\otimes_\CC {\rm
    End}_\CC\ (W_{\lambda '})\lra e_{\lambda '}\star \HH (G)\star
  e_{\lambda '}
$$
qui donn\'e par la m\^eme formule que celle de $\Upsilon_{\lambda
  '\omega}$.

 On a aussi un homomorphisme surjectif d'alg\`ebres :
$$
\Pi_{\omega}~: \ e_{\lambda '}\star \HH (G)\star e_{\lambda '}\lra
 \ e_{\lambda '\omega}\star \HH (G)\star e_{\lambda '\omega} 
$$
\noi donn\'e par
$$
\Pi_{\omega}(f)=\int_{F^{\times}}\omega (z)f(zg)\ d\mu_{F^\times}(z)\ .
$$

 On v\'erifie facilement le r\'esultat suivant.

\begin{lemma} Le diagramme suivant est commutatif :
$$
\begin{array}{ccccc}

  &                 & \Upsilon_{\lambda '} & & \\

   &   \HH (G,\lambda ' )\otimes {\rm End}_\CC \ (W_{\lambda '}) & \lra 
& e_{\lambda '}\star  \HH (G)\star e_{\lambda '} &   \\
P_{\omega}\otimes {\rm id} & \downarrow &  & \downarrow & \Pi_{\omega}
\\
  &  \HH (G,\lambda '\omega )\otimes {\rm End}_\CC\ (W_{\lambda '}) &
\lra & e_{\lambda '\omega}\star \HH (G)\star e_{\lambda '\omega} &\\
&  & \Upsilon_{\lambda '\omega} & &    
\end{array}
$$
\end{lemma}

\section{S\'erie discr\`ete et alg\`ebres de Hecke}

Soient $(J, \lambda )$ un type simple comme au {\S}1 
(dont nous gardons les notations) et $(J' ,\lambda ')$
son type modifi\'e. On note ${\mathcal R}(G)$ (resp. ${\mathcal R}(H)$) 
la cat\'egorie des repr\'esentations lisses de $G$ (resp. de $H$),
 ${\mathcal R}_\lambda (G)$ la sous-cat\'egorie pleine de ${\mathcal R}(G)$
 correspondant aux types $\lambda$ et $\lambda '$, et finalement
 ${\mathcal R}_{\rm Iw} (H)$ la sous-cat\'egorie pleine de ${\mathcal R}(H)$
 correspondant au type $(U(\Cfr ), {\mathbf 1}_{U(\Cfr )})$. 

\medskip

 Fixons un isomorphisme {\bf unitaire} d'alg\`ebres de Hecke
$$
\Psi\  : \ \HH (H,{\mathbf 1}_{U(\Cfr )}) \lra \HH (G,\lambda ')
$$
\noi comme dans le paragraphe {\S}3. Il induit une \'equivalence
 $\Psi_*$ entre les cat\'egories de modules \`a gauche correspondantes :
$$
\Psi_* \ : \  \HH (H,{\mathbf 1}_{U(\Cfr )})-{\rm Mod} 
\lra \HH (G,\lambda ')-{\rm Mod}\ .
$$
\noi Concr\`etement, si $M$ est un $\HH (H,{\mathbf 1}_{U(\Cfr )})$-module,    
   $\Psi_* (M)$ est le module b\^ati sur le m\^eme $\CC$-espace avec l'action :
$$
f.m = \Psi^{-1}(f).m\ , \ f\in \HH (G,\lambda ')\ , \ m\in M\ .
$$

Rappelons aussi les \'equivalences de cat\'egories induites par
 la {\it Th\'eorie des Types} (nous renvoyons le lecteur \`a [BKTypes]):
$$
\begin{array}{ccccc}
{\rm Mod}_\lambda & : & {\mathcal R}_{\lambda} (G) & \lra &
 \HH (G,\lambda ')-{\rm Mod} \\
                &   &  (\pi ,{\mathcal V}) & \mapsto &
 {\rm Hom}_{J'} (\lambda ' ,{\mathcal V})
\end{array}
$$

$$
\begin{array}{ccccc}
{\rm Mod}_{\rm Iw} & : & {\mathcal R}_{\rm Iw} (H) & \lra & 
\HH (H, {\mathbf 1}_{U( \Cfr )})-{\rm Mod} \\
                &   &  (\sigma ,{\mathcal W}) & 
\mapsto & {\rm Hom}_{U(\Cfr )} ({\mathbf 1}_{U(\Cfr ) } ,{\mathcal W})
 \simeq {\mathcal W}^{U(\Cfr )}
\end{array}
$$

 On note ${\mathcal E}_{\Psi}$ l'\'equivalence de cat\'egories rendant le diagramme 
suivant commutatif :
$$
\begin{array}{ccccc}
   &                           &  {\mathcal E}_\Psi &                       & \\
   & {\mathcal R}_{\rm Iw}(H )   & \lra              & {\mathcal
    R}_\lambda (G) &  \\
{\rm Mod}_{\rm Iw} & \downarrow &                    & \downarrow & 
{\rm Mod}_\lambda \\
   & \HH (H,{\mathbf 1}_{U(\Cfr )})-{\rm Mod} & \lra &
 \HH (G,\lambda ')- {\rm Mod} & \\
  &                                        & \Psi_* &
                            &
\end{array}
$$

 Rappelons qu'une repr\'esentation lisse irr\'eductible d'un 
groupe r\'eductif $p$-adique 
est dite de  {\it carr\'e int\'egrable} si son caract\`ere central
 est unitaire et si ses
coefficients sont de carr\'e int\'egrable modulo le  centre. 
Nous r\'eserverons 

\noi l'expression
 {\it repr\'esentation de la s\'erie discr\`ete} (ou encore {\it repr\'esentation
 essentiellement de carr\'e int\'egrable}) \`a une repr\'esentation obtenue
 d'une
 repr\'esentation 
de carr\'e int\'egrable par torsion par un caract\`ere lisse du groupe.
\medskip

 Nous aurons besoin du r\'esultat suivant d\^u \`a Bushnell et Kutzko. 

\begin{theorem} 
([BK] Theorem (7.7.1) et Gloss (7.7.2), page 257.) Soit $(\sigma ,W)$ 
une repr\'esentation 
irr\'eductible, objet de ${\mathcal R}_{\rm Iw}(H)$. Alors $(\sigma ,W)$
 est de carr\'e int\'egrable
si, et seulement si, ${\mathcal E}_\Psi (\sigma )$ est de carr\'e int\'egrable.
\end{theorem}

 Soit $B\subset H$ le sous-groupe de Borel des matrices triangulaires 
sup\'erieures. Notons 
$i_B^H$ et $I_B^H$ les foncteurs d'induction parabolique
respectivement  normalis\'e et non normalis\'e. Soit $\sigma$ une
repr\'esentation  de carr\'e int\'egrable de $H$, suppos\'ee 
Iwahori-sph\'erique (objet de ${\mathcal R}_{\rm Iw}(H)$). On sait 
qu'il existe un unique caract\`ere  unitaire non ramifi\'e de $K^{\times}$
 tel que $\sigma$ soit isomorphe \`a l'induite parabolique
$$
i_B^H \ (\vert \ \vert_K^{\frac{1-n}{2}}\chi \otimes \vert \ 
\vert_K^{\frac{3-n}{2}}\chi \otimes \cdots 
\otimes  \vert \ \vert_K^{\frac{n-1}{2}}\chi ) =
 I_B^H\ (\chi \otimes \cdots \otimes \chi ) =
\chi\circ{\rm det}\otimes I_B^H {\mathbf 1}_B\ .
$$

\noi La repr\'esentation de Steinberg de $H$, not\'ee ${\rm St}_H$ est 
par d\'efinition l'unique quotient  irr\'eductible de $I_B^H\ {\mathbf
  1}_B$.  Les autres repr\'esentations Iwahori-sph\'eriques et de carr\'e
 int\'egrable de $H$ sont donc, \`a isomorphisme pr\`es, les tordues
 $\chi\otimes {\rm St}_H$, o\^u $\chi$ parcourt les caract\`eres unitaires
 non ramifi\'es de $K^\times$. 

 Nous avons donc montr\'e le r\'esultat suivant.

\begin{lemma} Les repr\'esentations de ${\mathcal R}_\lambda (G)$ qui
  sont  irr\'eductibles et  de carr\'e int\'egrable sont, \`a isomorphisme
 pr\`es, les 
$$
\pi_\chi ={\mathcal E}_\Psi (\chi\otimes {\rm St}_H)
$$
\noi o\`u $\chi$ d\'ecrit les caract\`eres unitaires et non ramifi\'es de $K^\times$.
\end{lemma}

Soit $\chi$ un caractère unitaire et non ramifié de $K^\times$. 
 Un calcul immédiat, basé sur l'isomorphisme d'algèbres de Hecke
 de [BK](5.6.6), montre qu'il existe un isomorphisme unitaire d'algèbres de Hecke
$$
\Psi'~: \ \HH (H,{\mathbf 1}_{U(\Cfr )}) \lra \HH (G,\lambda ')
$$
tel que 
$$
{\mathcal E}_\Psi (\chi\otimes {\rm St}_H )\simeq {\mathcal E}_{\Psi '} ({\rm St}_H )
$$

\begin{corollary} Soit $\pi$ une représentation irréductible de carré intégrable, objet de 
 ${\mathcal R}_\lambda (G)$. Il existe alors un isomorphisme unitaire d'algèbres de Hecke
$$
\Psi~: \ \HH (H,{\mathbf 1}_{U(\Cfr )}) \lra \HH (G,\lambda ')
$$
tel que 
$$
\pi \simeq {\mathcal E}_{\Psi } ({\rm St}_H )
$$
\end{corollary}

\section{Le principe du transfert d'un pseudo-coefficient}

 Pour all\'eger les notations, nous notons $I=U(\Cfr )$. 

 Rappelons qu'on a une \'equivalence de cat\'egories :
$$
\HH (H, {\mathbf 1}_I)\text{-Mod} \lra \Rep_{\rm Iw}(H)
$$
\noi dont l'inverse est donn\'ee par
$$
\begin{array}{cccc}
{\rm Mod_I}~: & \Rep_{\rm Iw}(H) & \lra & \HH (H,{\mathbf 1}_I )\text{-Mod}\\
              & (\pi , {\mathcal V}) & \mapsto & {\mathcal V}^I \simeq
{\rm Hom}_I ({\mathbf 1}_I ,{\mathcal V})
\end{array}
$$

On fixe une repr\'esentation irr\'eductible $(\sigma_0 ,\WW_0 )\in \Rep_{\rm
  Iw}(H)$, de carr\'e int\'egrable et de caract\`ere central (unitaire)  fix\'e
$\omega_0$. On note $M_0 =\WW_0^I$ le $\HH (H,{\mathbf 1}_I)$-module
correspondant ; on sait qu'il est de dimension $1$. 
\medskip

 Soit $f_0$ un pseudo-coefficient de $(\sigma_0 ,\WW_0 )$. Par
 d\'efinition, $f_0$ est un \'el\'ement de $\HH (H,\omega_0 )$ qui v\'erifie
 la condition suivante. Pour toute repr\'esentation irr\'eductible lisse 
 temp\'er\'ee $(\sigma ,\WW )$ de $H$, de caract\`ere central $\omega_0$, on
 a 
$$
{\rm Tr}(\sigma (f_0 ), \WW )=
\left\{
\begin{array}{ll}
1  & \text{si }\sigma \simeq \sigma_0\\
0  & \text{si } \sigma \not\simeq \sigma_0
\end{array}
\right.
$$

 Nous
 faisons l'hypoth\`ese suivante :

\begin{hypothese} La fonction $f_0$ appartient \`a $\HH
  (H,\omega_0 {\mathbf 1}_I )$.
\end{hypothese}

 Soit $F_0\in \HH (H,{\mathbf 1}_I )$ une fonction telle que
 $P_{\omega_0}(F_0 ) =f_0$. 

 Notons que si $(\sigma ,\WW )$ est une repr\'esentation de $\Rep_{\rm
   Iw}(H)$ de caract\`ere central $\omega_0$, on a successivement :
$$
\begin{array}{ll}
{\rm Tr}(\sigma (f_0 ),\WW ) &  ={\rm Tr}(\sigma (F_0 ) ,\WW ) \text{, 
  par le lemme (4.1)}\\
                             &  ={\rm Tr} (\sigma (F_0 \star e_I ),
\WW ) \\
                             & = {\rm Tr} (\sigma (F_0 ) , e_I \star
\WW )\\
                             & = {\rm Tr}(F_0 ,M )
\end{array}
$$
\noi o\^u $M={\rm Mod}_I (\WW ) =\WW^I$. On a donc le 

\begin{lemma} Pour toute repr\'esentation temp\'er\'ee $(\sigma ,\WW )$ de
  $\Rep_{\rm Iw}(H)$, de caract\`ere central $\omega_0$, et de module
  ${\rm Mod}_I (\WW ) = M$, on a 
$$
{\rm Tr}(F_0 ,M) = 
\left\{
\begin{array}{ll}
1 & \text{si } \sigma\simeq \sigma_0\text{ c'est-\`a-dire si } M\simeq
M_0 \\
0 & \text{si } \sigma\not\simeq \sigma_0 \text{ c'est-\`a-dire si }
M\not\simeq M_0
\end{array}
\right.
$$
\end{lemma}

 Fixons un isomorphisme unitaire d'alg\`ebres de Hecke
$$
\Psi~: \ \HH (H,{\mathbf 1}_I )\lra \HH (G,\lambda ')\ .
$$
\noi Il induit une \'equivalence de cat\'egories :
$$
\Psi^*~: \HH (H, {\mathbf 1}_I )\text{-Mod}\lra \HH (G,\lambda
')\text{-Mod}
$$
\noi o\^u pour $M\in  \HH (H, {\mathbf 1}_I )\text{-Mod}$, $\Psi^*
(M)=M$ comme $\CC$-espace vectoriel, et o\^u la structure de module est
donn\'ee par :
$$
\varphi \star m =\Psi^{-1} (\varphi )\star m\ , \ m\in M\ ,
\ \varphi\in \HH (G,\lambda ')\ .
$$
\noi En particulier, pour tout $\varphi_0 \in \HH (H,{\mathbf 1}_I )$
et tout $M\in \HH (H, {\mathbf 1}_I )$-Mod, on a 
$$
{\rm Tr}(\varphi_0 , M) = {\rm Tr}(\Psi (\varphi_0 ),\Psi^* (M))\ .
$$

Le r\'esultat suivant est l'outil principal qui fait fonctionner notre
proc\'edure. Il se d\'eduit ais\'ement des travaux de Bushnell, Henniart
et Kutzko [BHK] sur le lien entre types et formule de Plancherel et sera
d\'emontr\'e dans l'annexe B. 

\begin{theorem} Soit $(\sigma ,\WW ) \in \Rep_{\rm Iw}(H)$ une
  repr\'esentation irr\'eductible de $\HH (H, {\mathbf 1}_I )$-module $M$
  et soit $(\pi ,\VV) = {\mathcal E}_\Psi (\sigma ,\WW )$ la
  repr\'esentation de $\Rep_{\lambda}(G)$ de $\HH (G,\lambda ')$-module
    associ\'e $\Psi^* (M)$. Alors la repr\'esentation $(\sigma  ,\WW )$ est
    temp\'er\'ee si, et seulement si, la repr\'esentation $(\pi ,\VV )$
    l'est.
\end{theorem}

 Notons $M_\lambda =\Psi^* (M_0 )$, et soit $(\pi_\lambda ,\VV_\lambda
 )={\mathcal E}_\Psi (\sigma_0 ,\WW_0 )$ la repr\'esentation de
 $\Rep_\lambda (G)$ correspondant au module $M_\lambda$. Par le
 Th\'eor\`eme (5.1), la repr\'esentation $(\pi_\lambda ,\VV_\lambda )$ est
 de carr\'e int\'egrale. Notons $\omega$ son caract\`ere central.

\begin{notation} Nous consid\'erons les fonctions suivantes :
\medskip

i) $F_\lambda = \Psi (F_0 )\in \HH (G,\lambda ')$,
\smallskip

ii) $f_\lambda = P_{\omega} (F_{\lambda})\in \HH (G , \lambda '\omega
)$,
\smallskip

iii) $\ds  \varphi_\lambda =\Upsilon_{\lambda '\omega}\big\{  \frac{1}{{\rm dim}\ (\lambda ')}
\ f_\lambda \otimes {\rm id}_{W_{\lambda }} \big\} ={\rm Tr}_{W_{\lambda '}}\circ f_{\lambda} \in e_{\lambda
  '\omega}\star \HH (G,\omega )\star e_{\lambda '\omega}$.
\end{notation}

\begin{proposition} La fonction  $\varphi_\lambda$ est
 un pseudo-coefficient de la repr\'esentation $(\pi_\lambda ,\VV_\lambda )$. 
\end{proposition}

\noi  {\it D\'emonstration} Soit $(\pi ,\VV )$ une repr\'esentation lisse irr\'eductible temp\'er\'ee, 
de caract\`ere central $\omega$,  du groupe $G$. Distinguons deux cas.
\medskip

\noi {\it Cas no 1}.  {\it On a} $(\pi ,\VV )\not\in \Rep_{\lambda}(G)$.
 En particulier $\pi\not\simeq \pi_\lambda$. Par d\'efinition 
$\pi (e_{\lambda '})\VV =\pi (e_{\lambda '\omega}).\VV =0$. On a donc
$$
{\rm Tr}\ (\pi (\varphi_\lambda ),\VV )={\rm Tr}\ (\pi (\varphi_\lambda \star e_{\lambda '\omega}),\VV )
 = {\rm Tr}\ (\pi (\varphi_\lambda ), \pi (e_{\lambda '\omega}).\VV )=0
$$

\noi {\it Cas no 2}.  {\it Supposons que} $(\pi ,\VV )\in \Rep_\lambda
(G)$.  Notons :
\medskip

 -- $M$ le $\HH (G,\lambda)$-module correspondant \`a $(\pi ,\VV )$ ;
\smallskip

 -- $(\sigma, \WW ) ={\mathcal E}_\Psi^{-1}  (\pi ,\VV )$ ;
\smallskip

 -- $N$ le $\HH (H,{\mathbf 1}_I )$-module correspondant \`a $(\sigma ,\WW )$. 
\medskip

 Par hypoth\`ese, la repr\'esentation $(\sigma ,\WW )$ est temp\'er\'ee.
\medskip

 Rappelons que les alg\`ebres $\HH (G,\lambda ')$ et $e_{\lambda '}
 \star \HH (G)\star e_{\lambda '}$ 
(resp. $(\HH , \lambda '\omega )$ et
 $e_{\lambda \omega}\star \HH (G,\omega )\star e_{\lambda '\omega })$)
 sont  \'equivalentes au sens de Morita, et que 
l'on a l'on via l'isomorphisme $\Upsilon_{\lambda '\omega}$ :
$$
\VV^{\lambda '}=\VV^{\lambda ' \omega}\simeq M\otimes W_{\lambda '}
$$

Nous pouvons \'ecrire successivement :

\begin{eqnarray}
{\rm Tr}\, (\varphi_\lambda ,\VV ) &=& {\rm Tr}\, (\varphi_\lambda
,\VV^{\lambda '})\\
 &=& {\rm Tr}\, \big( \Upsilon_{\lambda '\omega}\big\{ \frac{1}{{\rm dim}\,
(\lambda ')} f_\lambda \otimes {\rm id}_{W_{\lambda '}} \big\}
  ,\VV^{\lambda '}\big)\\
 &=& {\rm Tr}\, \big(  \frac{1}{{\rm dim}\,
(\lambda ')} f_\lambda \otimes {\rm id}_{W_{\lambda '}} ,M\otimes
    W_{\lambda '}\big)\\
 &=& {\rm Tr}\, (f_\lambda , M)\, \frac{1}{{\rm dim}\, (\lambda ')}
      {\rm Tr}({\rm id}_{W_{\lambda '}}, W_{\lambda '})\\
 &=& {\rm Tr}(f_\lambda ,M)\\
 &=& {\rm Tr}(F_\lambda ,M)\\
 &=& {\rm Tr}(\Psi (F_0 ),\Psi^* (N))\\
 &=& {\rm Tr}(F_0 ,N)
 \end{eqnarray}

Soit $\omega_0 ' = (\omega_0 )_{\vert F^\times}$ le caract\`ere trivial de
$F^\times$. Notons
$$
P_{\omega_0 '}~: \ \HH (H)\lra \HH (H,\omega_0 ')
$$
\noi le morphisme d'alg\`ebres donn\'e par
$$
P_{\omega_0 '}(f)(h)=\int_{F^\times} \omega_0 '(z)f(zh)\,
d\mu_{F^\times} (z)\ , _ f\in \HH (f) , \ h\in H\ ,
$$
\noi o\^u $\mu_{F^\times}$ est la mesure de Haar sur $F^\times$ donnant
la volume $1$ \`a $\ofr_F^\times$.

 Notons $f_0 ' = P_{\omega_0 '}(F_0 )$. Par le lemme (7.7.6) de [BK],
  le caract\`ere central de $\sigma$ v\'erifie :
 $(\omega_{\sigma})_{\vert F^\times} = \omega_0 '$, de sorte que
 $\sigma$ peut se voir comme une repr\'esentation de $H/F^\times$. En
 utilisant le crit\`ere $L^{2+\epsilon}$, on voit que $\sigma$ est une
 repr\'esentation temp\'er\'ee de $H/F^\times$. 

 Nous faisons la {\it seconde hypoth\`ese} suivante :

\begin{hypothese}  Il existe une constante non nulle $c\in \CC$, telle
  que, vue comme repr\'esentation de $H/F^\times$,
  $(\sigma_0 , \WW_0 )$ admet la fonction $c. f_0 '$ comme
  pseudo-coefficient.
\end{hypothese}

 On obtient ainsi :

\begin{eqnarray}
{\rm Tr}\, (\varphi_\lambda ,\VV ) &=& {\rm Tr}\, (F_0 ,\WW_\sigma )\\
 &=& \frac{1}{c}.{\rm Tr}\, (cf_0 ' , \WW_{\sigma}) \ .
\end{eqnarray}

\noi On en d\'eduit que si $\sigma\not\simeq \sigma_0$ (i.e. si
$\pi\not\simeq\pi_\lambda$), on a ${\rm Tr}\, (\varphi_\lambda ,\VV
)=0$.

D'un autre c\^ot\'e, si $\sigma\simeq \sigma_0$ (i.e. $\pi \simeq
\pi_\lambda$), on a alors 

\begin{eqnarray}
{\rm Tr}\, (\varphi_\lambda ,\VV ) &=& {\rm Tr}\, (F_0
,\WW_{\sigma_0})\\
 &=& {\rm Tr}(f_0 ,\WW_{\sigma_0})\\
 &=& 1
\end{eqnarray}

\noi ce qui termine notre d\'emonstration.

\section{Le pseudo-coefficient de Kotwittz}

 On suit ici la section 2. de [Kott], o\^u Kottwitz d\'efinit des fonctions
 d'{\it Euler-Poincar\'e} $f_{\rm EP}$ pour tout groupe r\'eductif
 connexe \`a centre  anisotrope. 
\medskip

 Soit donc ${\mathbb L}$ un groupe r\'eductif connexe de centre
  anisotrope. On note $ L$ le
 groupe de ses  points $F$-rationnels.   Fixons une mesure
 de Haar $\mu_{L}$ sur $L$ et
 notons :
\medskip

 -- $X_{L}$ l'immeuble de Bruhat-Tits de $\mathbb L$ sur $F$ ;
\smallskip

 -- $\Sigma$ un syst\`eme de repr\'esentants des $L$-classes de
 conjugaison de simplexes de $X_{L}$ ;
\smallskip

 -- $d_{L}$ le $F$-rang de $\mathbb L$ ;
\smallskip

 -- $d_\sigma$ la dimension d'un simplexe $\sigma$ de $X_{L}$ ;
\smallskip

 -- ${L}_\sigma$ le stabilisateur d'un simplexe $\sigma$ dans
 $L$ ;
\smallskip

 -- ${\rm sgn}_\sigma (x)$ la signature de la permutation des sommets
 d'un simplexe $\sigma$ de $X_{L}$ induite par l'action de $x\in
 {L}_\sigma$.

\smallskip

 -- ${\mathbf 1}_{U}$ la fonction caract\'eristique d'une partie $U$ de
 $L$ ;
\medskip

 Kottwitz d\'efinit une fonction d'Euler-Poincar\'e par la formule :
\medskip

$$
f_{\rm EP}^{\Sigma , L} = \sum_{\sigma \in \Sigma}
(-1)^{d_\sigma } \frac{1}{\mu_{L}({L}_\sigma)} {\mathbf
  1}_{{L}_\sigma} \, {\rm sgn}_\sigma\ .
$$

On a alors le r\'esultat fondamental suivant.

\begin{theorem} (Kottwitz-Casselman, [Kott], Thm 2', page 637.) La
  fonction $f_{\rm Kottwitz}^{\Sigma ,L} := (-1)^{d_L -1} f_{\rm
    EP}^{\Sigma ,L}$ est un pseudo-coefficient de la repr\'esentation de
  Steinberg de $L$.
\end{theorem}

 Nous allons sp\'ecialiser ce r\'esultat \`a deux groupes $\mathbb L$
 particulier. Dans le premier cas ,  c'est $K$ qui joue le r\^ole du corps
 de base $F$.
\medskip

 Commen\c cons par le $K$-groupe r\'eductif  ${\mathbb L}= {\rm PGL}(e)$. Ici~:
\medskip

 -- la mesure de Haar $\mu_L$ est prise comme \'etant $\mu_H
 /\mu_{K^\times}$, o\^u $\mu_H$ est la mesure de Haar sur $H={\rm
   GL}(e,K)$ qui donne le
 volume $1$ \`a un sous-groupe d'Iwahori et $\mu_{K^\times}$ la mesure
 de Haar sur $K^\times$ qui donne la mesure $1$ \`a $\ofr_K^\times$.
\smallskip

 -- $\Sigma$ est un syst\`eme de repr\'esentants des $H$-orbites de
 simplexes construit comme dans la section 2. 
\medskip

 En particulier, on note~:
\medskip

 -- $\Theta$ un syst\`eme de repr\'esentants des orbites de $\Pi$ dans
 l'ensemble des parties de $S$ ;
\smallskip

 -- ${\mathcal K}_T$ le normalisateur dans $H$ du parahorique
 $P_T=U(\Cfr )\langle T\rangle U(\Cfr )$, et ${\bar {\mathcal K}}_T$
 son image dans ${\rm PGL}(e,K)$, et ${\mathbf 1}_{{\mathcal K}_T}$ sa
 fonction caract\'eristique ;
\smallskip

 -- $\sigma_T$ l'unique simplexe de l'immeuble fix\'e par ${\mathcal
   K}_T$, et $d_T$ sa dimension, i.e. $d_T =e-1 -\vert T\vert$. 
\smallskip

 -- ${\rm sgn}_T ={\rm sgn}_{\sigma_T}$.
\medskip

 On a alors :
$$
f_{\rm EP}^{\Sigma ,{\mathbb L}} =f_{\rm EP}^{\Theta ,H} := 
\sum_{T\in \Theta}(-1)^{d_T} \frac{1}{\mu_{H/K^\times} ({\bar {\mathcal
      K}}_T )} {\mathbf 1}_{{\bar {\mathcal K}}_T } {\rm sgn}_T\ .
$$

Regardons $f_{\rm EP}^{\Theta ,H}$ comme une fonction de l'alg\`ebre de
Hecke $\HH (H,\omega_0 )$.

\begin{corollary} La fonction $f_{\rm Kottwitz}^{\Theta ,H} :=
  (-1)^{e-1} f_{\rm EP}^{\Theta ,H}$ est un pseudo-coefficient de la
  repr\'esentation de Steinberg de $H$.
\end{corollary}

 Dans [La] {\S}5, G\'erard Laumon propose une variante de ce
 pseudo-coefficient que nous allons rappeler. Soit $f_0$ la fonction
 de $\HH (H,\omega_0 )$ donn\'ee par
$$
f_0 = (-1)^{e-1} \sum_{T\subset S} (-1)^{d_T} \frac{{\rm sgn}_T  . {\mathbf
    1}_{{\mathcal K}_T} }{(d_T +1){\rm
    vol}(P_T , dh)} \ .
$$

Nous avons alors le

\begin{proposition} ([La], Lemma (5.2.2), page 135) 

\noi (i) La fonction $f_0$ est la moyenne des fonctions $f_{\rm
  Kottwitz}^\Theta$, lorsque $\Theta$ d\'ecrit tous les ensembles de
  repr\'esentants possibles des orbites de $\Pi$ dans l'ensemble des
  parties de $S$.
\smallskip

 \noi (ii) La fonction $f_0$ est un pseudo-coefficient de la
 repr\'esentation de Steinberg de $H$.
\end{proposition}

Bien s\^ur le point (ii) d\'ecoule de (i) de fa\c con imm\'ediate. 
\medskip

 Pour $T\subset S$, le normalisateur de $P_T$ dans $H$ peut s'\'ecrire
 $z_T^{\mathbb Z} P_T$, o\^u $z_T$ est une certaine puissance $\Pi^u$ de
 $\Pi_T$, $u$ entier positif divisant $e$. 

\begin{notation} On pose $\epsilon_T ={\rm sgn}_T (z_T )$ et on note $n_T$
  l'unique entier positif tel que $z_T^{n_T} =\varpi_K$·
\end{notation}

\begin{lemma} Pour $h\in H$, on  a l'\'egalit\'e
$$
f_0 (h) =(-1)^{e-1}\sum_{T\subset S}\frac{(-1)^{d_T}}{(d_T +1)\mu_H (P_T )}
\sum_{l=0,...,n_T -1} \epsilon_T^{l} \sum_{w\in \langle T\rangle } \sum_{k\in \ZZ} 
f_{w}^0 (z_T^{-l}\varpi_K^k h)\ .
$$
\end{lemma}

\noi {\it D\'emonstration} Cette formule d\'ecoule facilement des
\'ecritures :
$$
{\mathcal K}_T =\coprod_{k\in \ZZ} z_T^k P_T
$$
\noi et 
$$
{\mathbf 1}_{P_T} = \sum_{w\in \langle T\rangle} {\mathbf 1}_{IwI} =
 \sum_{w\in \langle T\rangle} f_w^0\ .
$$

 Pour la suite nous aurons besoin d'une fonction $F_0 \in \HH (H)$,
 telle qu'avec les notations de la section 4, on ait
 $P_{\omega_{0}}(F_0 ) = f_0$. Il est clair que l'on peut se donner
 une telle $F_0$ par la formule suivante :

\begin{equation}
F_0 (h) =  \frac{(-1)^{e-1}}{e'}\sum_{T\subset S}
\frac{(-1)^{d_T}}{(d_T +1)\mu_H (P_T )}
\sum_{w\in \langle T\rangle} \sum_{l=0,...,e'n_T -1} \epsilon_T^{l}
  f_{w}^0 (z_T^{-l} h)\ , \ h\in H\ ,
\end{equation}

\noi o\^u $e' =e(E/F)=e(K/F)$.
\bigskip

 Le second groupe auquel nous allons sp\'ecialiser le r\'esultat du
 Th\'eor\`eme (7.1) est $H' = {\rm GL}(e,K)/F^\times$, o\^u $F^\times$ est
 naturellement vu comme un sous-groupe du centre $K^\times$ de $H$.
On consid\`ere ici $H'$ comme le groupe des $F$-points rationnels du
$F$-groupe r\'eductif connexe ${\rm Res}_{K/F}\big(  {\rm GL}(e)/K \big)
/{\rm GL}(1)/F$.
Le centre $K^\times /F^\times$ de $H'$ \'etant compact, nous sommes sous
les hypoth\`eses du th\'eor\`eme (7.1). 
\medskip

 On munit $H'$ de la mesure de Haar obtenue en quotientant $\mu_H$ par
 $\mu_{F^\times}$, la mesure de Haar sur $F^\times$ donnant le volume
 $1$ \`a $\ofr_F^\times$. Le Th\'eor\`eme (7.1) fournit des
 pseudo-coefficients $f_{\rm Kottwitz}^{\Sigma ,H'}$ pour la
 repr\'esentation de Steinberg de $H'$. Nous allons comparer ces
 pseudo-coefficients \`a ceux de la repr\'esentation de Steinberg de ${\bar
 H}={\rm PGL}(e,K)$. 
 
\medskip

\begin{lemma} Les groupes $\bar H$ et $H'$ ont le m\^eme immeuble. 
\end{lemma}

\noi {\it D\'emonstration}. En effet, ${\bar H}={\rm PGL}(e,K)$ est obtenu 
en quotientant $H'$ par son   centre. Les immeubles  de ces deux groupes
s'identifient donc naturellement. 

\begin{lemma} Notons $p$~: $H' \lra {\bar H}$ la projection naturelle et
  $\St_{\bar H}$, $\St_{H'}$ les repr\'esentations de Steinberg de $\bar
  H$ et $H'$   respectivement. On peut les r\'ealiser dans le m\^eme espace o\`u elle
  sont reli\'ees par la relation :
$$
\St_H (h) = \St_{H'}(p(h))\ , \ h\in {\bar H} \ .
$$
\end{lemma}

\noi {\it D\'emonstration}. Par [Borel-Serre], $\St_{\bar H}$
(resp. $\St_{H'}$) se r\'ealise comme la repr\'esentation naturelle de
$\bar H$ (resp. de $H'$) dans l'espace de cohomologie \`a support compact et
coefficients complexes $H_{c}^{e-1} (X,\CC )$, o\^u $X$ d\'esigne
l'immeuble. Le r\'esultat est donc \'evident. [Il y a s\^urement plus simple
  !!]
\medskip

 Fixons un syst\`eme de repr\'esentants $\Sigma$ des $\bar H$-classes de
 conjugaison de simplexes dans l'immeuble $X$ de $\bar H$ (et de $H'$
 !).  Par le lemme (7.6), c'est aussi un syst\`eme de repr\'esentants des
 $H'$-classes de  conjugaison de l'immeuble de $H'$.

On peut voir la fonction d'Euler-Poincar\'e $f_{\rm EP}^{\Sigma , H'}$
comme une fonction de $\HH (H,{\mathbf 1}_{F^\times})$. Consid\'erons le
morphisme d'alg\`ebres 
$$
P'~: \ \HH (H, {\mathbf 1}_{F^\times})\lra \HH (H, {\mathbf
  1}_{K^\times})
$$
\noi donn\'e par
$$
P' (f)(h) = \int_{K^\times /F^\times} f(zh)\, d\mu_{K^\times
    /F^\times}(z) \ , \  h\in H\ .
$$
\noi Ici $\mu_{{K^\times /F^\times}}$ est la mesure de Haar sur le
  centre $K^\times /F^\times$ quotient de la mesure de Haar
  $\mu_{K^\times}$ par $\mu_{F^\times}$. 
\bigskip

 Par un calcul immédiat laissé au lecteur, nous avons le résultat suivant :

\begin{lemma} Les fonctions $f_{\rm EP}^{\Sigma , H'}$ et
  $f_{EP}^{\Sigma ,{\bar H}}$ sont reli\'ees par la relation :
$$
f_{EP}^{\Sigma ,{\bar H}} = c.P' \big(  f_{\rm EP}^{\Sigma , H'}\big)
\ ,
$$
\noi o\^u $c\in \CC$ est une constante non nulle.
\end{lemma}

 Pour un syst\`eme de repr\'esentants $\Sigma$ construit comme dans la
 section 2, on obtient :
$$
f_{\rm EP}^{\Sigma ,H'}=f_{\rm EP}^{\Theta ,H'} := \sum_{T\in \Theta}
(-1)^{d_T}\frac{1}{\mu_{H'}({\mathcal K}_T ')} {\mathbf 1}_{{\mathcal K}_T '}
  {\rm sgn}_T
$$
\noi o\^u ${\mathcal  K}_T '$ est l'image de ${\mathcal K}_T$ dans $H'$.

 Suivant Laumon, nous d\'efinissons une nouvelle fonction de $\HH (H ,
 {\mathbf 1}_{F^\times})$ par
 la formule :
$$
f_0 ' = (-1)^{e-1}\sum_{T\in S} (-1)^{d_T} \frac{{\rm sgn}_T . {\mathbf
    1}_{{\mathcal K}_T}}{(d_T +1 )\mu_{H}(P_T )}\ . 
$$

\begin{proposition} (i) Pour tout $\Theta$, $f_{\rm Kottwitz}^{\Theta
    ,H'} := (-1)^{e-1} f_{\rm EP}^{\Theta ,H'}$ est un pseudo
  coefficient de la repr\'esentation de Steinberg de $H'$.
\smallskip

\noi (ii) La fonction $f_0 '$ est la moyenne des fonctions $f_{\rm
  Kottwitz}^{\Theta ,H'}$ lorsque $\Theta$ d\'ecrit tous les ensembles
possible de repr\'esentants des orbites de $\Pi$ dans l'ensemble des
parties de $S$.
\smallskip

\noi (iii) En particulier, $f_{0}'$ est un pseudo-coefficient de la
repr\'esentation  de Steinberg de $H'$.
\end{proposition}

 Notons 
$$
P_{\omega_0}'~: \ \HH (H)\lra \HH (H' , {\mathbf 1}_{F^\times})
$$
\noi la projection naturelle donn\'ee par 
$$
P_{\omega_0}' (f)(g)= \int_{F^\times} f(zg)d\mu_{F^\times} (z)\ , \ f\in
\HH (H)\ .
$$

\begin{lemma} On a $P_{\omega_0}' (F_0 )=c. f_0 '$, pour une constante
  non nulle $c\in \CC$. En particulier si
  $(\sigma ,\WW)$ est une repr\'esentation lisse admissible de $H$ telle
  que 
$$
\sigma (zh) = \sigma (h)\ , \ h\in H , \ z\in F^{\times}\ ,
$$
on a alors

$$
{\rm Tr}\, (F_0 , \WW )= c . {\rm Tr}\, (f_0 ' , \WW )\ .
$$
\end{lemma}

\section{Transfert du pseudo-coefficient de Kotwittz}

Nous utilisons les notations des deux sections pr\'ec\'edentes.
Fixons une représentation $\pi =\pi_\lambda$ irréductible de carré intégrable, objet de 
${\mathcal R}_{\lambda}(G)$. On sait par le corollaire (5.3) qu'il existe un isomorphisme 
unitaire d'algèbres de Hecke 
$$
\Psi~: \ \HH (H,{\mathbf 1}_{U(\Cfr )}) \lra \HH (G,\lambda ')
$$
\noi tel que $\pi_\lambda$ s'écrive ${\mathcal E}_\Psi ({\rm St_H})$.  Le but de
cette section est de calculer le pseudo-coefficient $\varphi_\lambda$
de la proposition (6.5) obtenu par transfert du pseudo-coefficient de 
Kottwitz (plus exactement sa version modifiée de Laumon), au travers des isomorphismes
d'algèbres de Hecke.  Plus exactement nous allons d\'eterminer  la
  restriction de $\varphi_{\lambda}$ \`a $G^0$, o\`u  
$$
G^0 =\lbrace g\in G\ ; \ {\rm val}_F ({\rm det}\, (g))=0\rbrace\ .
$$ 

 Rappelons que

\begin{equation}
F_0  = \frac{(-1)^{e-1}}{e'}\sum_{T\subset S}\frac{(-1)^{d_T}}{(d_T
  +1)\mu_H (P_T)} \sum_{w\in \langle T\rangle}\sum_{l=0,...,e'n_T
  -1}\epsilon_T^l f_{z_{T}^{l} w}
\end{equation}

Le pseudo-coefficient $\varphi_\lambda$ se calcule \`a partir de $F_0$ par les trois
\'etapes :
\medskip

\noi (1) $F_\lambda =\Psi (F_0 )\in \HH (G,\lambda ')$,
\smallskip

\noi (2) $f_\lambda = P_{\omega} (F_\lambda )\in \HH (G,\lambda '\omega
)$,
\smallskip

\noi (3) $\varphi_\lambda = {\rm Tr}_{W_{\lambda '}}\circ f_\lambda \in
e_{\lambda '\omega}\star \HH (G,\omega )\star e_{\lambda '\omega}$.
\medskip

\noi o\^u $\omega$ d\'esigne le cacact\`ere central de $\pi_\lambda$.
\medskip

 Pour $g\in G$, on a alors :

\begin{equation}
\varphi_\lambda (g)= 
\frac{(-1)^{e-1}}{e'}\sum_{T\subset S}\frac{(-1)^{d_T}}{(d_T
  +1)\mu_H (P_T)} \sum_{w\in \langle T\rangle}\sum_{l=0,...,n_T e'
  -1}\sum_{k\in \ZZ}\omega (\varpi_F )^{-k}
\epsilon_T^l \Psi_{z_T^l  w}^k (g)
\end{equation}

\noi o\^u $\Psi_{z_T^l  w}^k$ est la fonction sur $G$ d\'efinie par 

\begin{equation}
 \Psi_{z_T^l  w}^k(g)=
{\rm  Tr}_{W_{\check{\lambda}'}}\circ \Psi (f_{z_T^l
  w})(\varpi_F^{-k}g)\ . 
\end{equation}

 Notons que la fonction $\Psi (f_{z_T^l w})$ a son support dans
 $J'z_T^l wJ'$, avec $J'\subset G^0$. Il s'ensuit que si $g\in G^0$,
 le terme  $\Psi (f_{z_T^l  w})(\varpi_F^{-k}g)$ est nul si $k\not= 0$
 ou $l\not= 0$. On a alors l'expression simplifi\'ee :
$$
\varphi_\lambda (g)=\frac{(-1)^{(e-1)}}{e'}\sum_{T\subset
  S}\frac{(-1)^{d_T}}{(d_T +1)\mu_H (P_T )}\sum_{w\in \langle
  T\rangle}\Psi_w^0 (g) , \ g\in G^0\ .
$$

 Pour expliciter le pseudo-coefficient, il nous reste donc \`a
 d\'eterminer ${\rm  Tr}_{W_{\check{\lambda}'}}\circ \Psi (f_{w}^0 )$,
 pour $w\in W^0$. Rappelons   que $\Psi (f_{w}^0)$ 
est l'\'el\'ement de $\HH
(G,\lambda ')$, de support $J' w J'$, donn\'e par
$$
\Psi (f_w^0 )(g)=\frac{1}{\vert \baP\vert}{\tilde \kappa}_M (g)\otimes
     [ {\tilde \sigma}({\bar j}_1)\circ T_w \circ {\tilde
         \sigma}({\bar j}_2)]\ ,
$$
\noi o\`u $g=j_1 wj_2$, $j_1$, $j_2\in J'$. Ici, pour $j\in J'$, $\bar
j$ d\'esigne l'image de $j$ dans le quotient $J' /J_M^1 =U(\Bfr )J_M^1
/J_M^1 \simeq \baP$. On a obtenu le r\'esultat suivant :

\begin{lemma} Pour $w\in W^0$, ${\rm Tr}_{W_{\check{W}'}}\circ \Psi
  (f_w^0 )$ est la fonction \`a support dans $J' wJ' \subset J_M$ donn\'ee
  par
\begin{equation}
{\rm Tr}_{W_{\check{W}'}}\circ \Psi
  (f_w^0 )(g)={\rm Tr}\, \check{\kappa}_M (g) . {\rm Tr}\, {\bar f}_w
({\bar g})
\end{equation}

\noi o\^u $\bar g$ d\'esigne l'image de $g$ dans $J_M /J_M^1$ et ${\bar
  f}_w$ la fonction introduite \`a la section 3.
\end{lemma}

 En particulier le support de $\varphi_{\vert G^0}$ est contenu dans
 la r\'eunion des supports des $\Psi (f_w)$, $w\in W^0$, c'est-\`a-dire :
$$
\bigcup_{w\in W^0} J^1_M  U(\Bfr )w U(\Bfr )J_M^1 \subset J_M^1
U(\Bfr_M )J_M^1 =J_M\ .
$$
\noi Nous avons donc d\'emontr\'e la formule suivante.

\begin{proposition} Soit $g\in G^0$. On a
$$
\varphi_\lambda (g)=\left\{
\begin{array}{ll}
0 & \text{si } g\not\in J_M\\
\frac{(-1)^{(e-1)}}{e'}\sum_{T\subset
  S}\frac{(-1)^{d_T}}{(d_T +1)\mu_H (P_T )}\sum_{w\in \langle
  T\rangle} {\rm Tr}\, {\tilde \kappa}_M (g).{\rm Tr}\, \bar{f}_w
(\bar{g}) & \text{si } g\in J_M\ .
\end{array}\right.
$$
\end{proposition}

\section{Calculs de traces dans les groupes finis}

 Les notations de cette section sont ind\'ependantes des pr\'ec\'edentes.
\bigskip

 On consid\`ere un groupe fini $G$ et une repr\'esentation complexe
 irr\'eductible $(\pi ,V)$ de $G$.

\begin{lemma} Supposons $V$ muni d'un produit scalaire
  $\langle\ , \ \rangle$ $G$-invariant. Soient $v\in V$ un
 vecteur norm\'e et $T\in
  {\rm End}_{\CC}(V)$. On a alors :
$$
{\rm Tr}(T)=\frac{{\rm dim}\ \pi}{\vert G\vert} \sum_{x\in G}\langle v,\pi
(x)\circ T\circ \pi (x^{-1}).v \rangle \ .
$$
\end{lemma}

\noi {\it D\'emonstration.} L'endomorphisme
$$
\frac{1}{\vert G\vert}\sum_{x\in G} \pi (x)\circ T\circ \pi (x^{-1})
$$
est $G$-\'equivariant et par le lemme de Schur, il est de la forme
$\lambda (T){\rm id}_V$, pour un $\lambda (T)\in \CC$. L'application
obtenue $\lambda$~: ${\rm End}_{\CC}(V)\lra \CC$ est une forme
lin\'eaire. On peut la voir comme une forme bilin\'eaire sur $V\otimes
{\tilde V}$, qui est $G$-invariante pour l'action de $G$ donn\'ee par
$\pi \otimes {\tilde \pi}$. Une seconde application du lemme de Schur
montre que $\lambda$ est proportionnelle \`a l'application trace : il
existe $k\in \CC$, tel que $\lambda =k. {\rm Tr}$. On a $\lambda ({\rm
  id}_V )=1$ et $k.{\rm Tr}({\rm id}_V )=k.{\rm dim}\ \pi$, de sorte
que $k=1/{\rm dim}\ \pi$.  Il vient donc
$$
\frac{{\rm dim}\ \pi}{\vert G\vert} \sum_{x\in G}(v,\pi
(x)\circ T\circ \pi (x^{-1}).v ) = \frac{1}{\vert G\vert}{\rm
  Tr}(T)\sum_{x\in G} (v,v) ={\rm Tr} (T)\ .
$$

 Notons comme cons\'equence, le r\'esultat bien connu suivant.

\begin{corollary}
 Soit $v\in V$ un vecteur de norme $1$
  pour un produit scalaire $\langle \ , \ \rangle$ $G$-invariant sur $V$. Alors  en
  notant $f$ le coefficient de $\pi$ donn\'e par $g\mapsto \langle v,\pi
  (g).v\rangle $, on a :
$$
{\rm tr}\pi(\gamma )=\frac{{\rm dim}\ \pi}{\vert G\vert}\sum_{x\in G}
f(x\gamma x^{-1})\ .
$$
\end{corollary}

\noi {\it D\'emonstration}. Appliquer le lemme \`a $T=\pi (x)$.
\bigskip

Donnons-nous \`a pr\'esent un sous-groupe $H$ de $G$ et $(\sigma ,V_H )$
une repr\'esentation irr\'eductible de $H$. Soit $(\pi ,V_G )$ la
repr\'esentation induite. On identifie comme d'habitude l'alg\`ebre
d'entrelacement ${\rm
  End}_G (V_G )$ \`a l'alg\`ebre de Hecke $\HH (G, {\tilde \sigma})$ des
fonctions $f$~: $G\lra {\rm End}_\CC (V_H )$ telle que $f(h_1 gh_2 )=\sigma
(h_1 )f(g)\sigma (h_2 )$, $g\in G$, $h_1$, $h_2 \in H$.

 Concr\`etement une fonction $\varphi\in \HH (G,{\tilde \sigma})$
 s'identifie \`a l'op\'erateur de convolution  $\varphi\in {\rm End}_G
 (V_G )$ donn\'e par $\varphi (f)=\varphi\star f$, o\`u
$$
\varphi \star f (g) =\sum_{x\in G} \varphi (x)[f(x^{-1}g)] =\sum_{x\in
  G}\varphi (gx^{-1})[f(x)]\ , \ f\in V_G\ , \ g\in G\ .
$$

 Fixons un idempotent $e$ de $\HH (G,{\tilde \sigma} )$ et notons
   $(\pi_e ,V_e )$ la sous-repr\'esentation de $(\pi ,V_G )$ donn\'ee par
$$
V_e =e\star {\rm Ind}_H^G \ \sigma = e\star V_G \ .
$$

Nous faisons les  hypoth\`eses suivantes.

\begin{hypothese} La repr\'esentation $(\pi_e ,V_e )$ est
  irr\'eductible.
\end{hypothese}

\begin{hypothese}
 Il existe un produit scalaire $H$-invariant
  $\langle \ , \ \rangle_H$ sur $V_H$  tel que,
  pour tout $x\in G$, l'adjoint $e^* (x)$ de $e(x)$ relativement \`a
  $\langle\ , \  \rangle_H$ est $e(x^{-1})$.
\end{hypothese}

\begin{hypothese} On a $e(1)=\lambda_1 {\rm id}_{V_H}$ pour
  une constante r\'eelle $\lambda_1 >0$.
\end{hypothese}

 Pour la suite nous fixons un tel produit scalaire.
\medskip

Nous nous
 proposons de d\'emontrer le r\'esultat suivant.

\begin{proposition} La trace de la repr\'esentation $\pi_e$ en
  un \'el\'ement $\gamma \in G$ est donn\'ee par la formule :
$$
{\rm Tr}\ \pi_e (\gamma )=\frac{1}{\lambda_1}\ \frac{{\rm
    dim}\ \pi_e}{{\rm dim}\ \sigma}\  \frac{\vert H\vert}{\vert G\vert}
\sum_{x\in H\backslash G} [{\rm Tr}\circ e](xyx^{-1})\ .
$$
\end{proposition}

\noi {\it D\'emonstration} Fixons un vecteur $v\in V_H$ tel que $\Vert
v\Vert_H =1$ et soit $f_v\in {\rm Ind}_H^G\ \sigma$ la fonction de
support $H$, donn\'ee par
$$
f_v (h)=\sigma (h).v\ , \ h\in H\ .
$$
\noi La fonction $g_e =e\star f_v$ est alors dans $V_e$. Pour $y\in
G$, on a
$$
g_e (y)=\sum_{x\in G} e(yx^{-1} )f_v (x)=\sum_{x\in H}e(y)\sigma
(x^{-1}) \sigma (x) f_v (1) = \vert H\vert . e(y)(v)\ .
$$

On d\'efinit un produit scalaire $G$-invariant sur $V_G$ par  la formule
$$
\langle f,g\rangle_G  =\sum_{x\in G}\langle f(x) ,g(x)\rangle_H\ ,
\ f,g\in V_G
$$
\noi On a
$$
\Vert g_e\vert_G^2 =\vert H\vert^2 \sum_{x\in G}\langle
e(x)(v),e(x)(v)\rangle_H
$$
$$
= \vert H\vert^2 \sum_{x\in G}\langle v,e(x^{-1})e(x) v\rangle_H
= \vert H\vert^2 \langle v,(\sum_{x\in G}e(x^{-1})e(x))(v)\rangle_H
$$
$$
= \vert H\vert^2 \langle v,e(1)(v)\rangle = \vert H\vert^2 \lambda_1\ .
$$
\noi La fonction $\displaystyle f_e =\frac{1}{\sqrt{\lambda_1}}g_e$, donn\'ee par
$\displaystyle f_e (y)=\frac{1}{\sqrt{\lambda_1}} e(y)(v)$ est donc un
vecteur norm\'e de $V_e$.

 Soit $\gamma\in G$. D'apr\`es le corollaire 9.2, on a
\bigskip

\begin{tabular}{rcl}
${\rm tr}\ \pi_e (\gamma )$ & $=$ & $\ds \frac{{\rm dim}\ \pi_e}{\vert G\vert}
  \sum_{x\in G}\langle f_e ,\pi_e (x\gamma x^{-1}) f_e\rangle_G$\\
      & $=$  &  $\ds \frac{{\rm dim}\ \pi_e}{\vert G\vert} \sum_{x\in G}
  \sum_{y\in G} \langle f_e (y),f_e (yx\gamma x^{-1})\rangle_H$\\
   & $=$ &  $\ds \frac{{\rm dim}\ \pi_e}{\lambda_1 . \vert G\vert} \sum_{x\in G}
  \sum_{y\in G} \langle e(y)(v),e(yx\gamma x^{-1})(v)\rangle_H$\\
  & $=$ &  $\ds \frac{{\rm dim}\ \pi_e}{\lambda_1 . \vert G\vert} \sum_{x\in
    G} \langle v,\big\{ \sum_{y\in G}e(y^ {-1})e(yx\gamma x^{-1})\big\}
  (v)\rangle_H$\\
  & $=$ &$  \ds \frac{{\rm dim}\ \pi_e}{\lambda_1 . \vert G\vert} \sum_{x\in
    G} \langle v,e\star e (x\gamma x^{-1})(v)\rangle_H$\\
  & $=$ &  $\ds \frac{{\rm dim}\ \pi_e}{\lambda_1 . \vert G\vert} \sum_{x\in
    G} \langle v,e (x\gamma x^{-1})(v)\rangle_H$
\end{tabular}
\bigskip

 \noi Mais d'apr\`es le lemme 9.1, on a pour tout $x\in G$ :
\begin{eqnarray}
\sum_{h\in H}\langle v , e(hx\gamma x^{-1} h^{-1}) v\rangle_H & = &
\sum_{h\in H}\langle v , \sigma (h) e(x\gamma x^{-1}) \sigma (h^{-1})
v\rangle_H \\
 & = & \frac{\vert H\vert}{{\rm dim}\ \sigma}{\rm Tr}\ e(x\gamma x^{-1})\ .
\end{eqnarray}
\noi On en d\'eduit donc bien :

$$
{\rm Tr}\ \pi_e (\gamma )  \frac{1}{\lambda_1}\ \frac{{\rm
    dim}\ \pi_e}{{\rm dim}\ \sigma}\  \frac{\vert H\vert}{\vert G\vert}
\sum_{x\in H\backslash G} [{\rm Tr}\circ e](xyx^{-1})\ .
$$

\section{Caract\`eres des repr\'esentations triviales g\'en\'eralis\'ees}

 Pour des preuves ou des r\'ef\'erences, nous renvoyons le lecteur \`a
 l'article [SZ] de Silberger et Zink.
\bigskip

Le groupe $\bar G$ poss\'ede deux repr\'esentations irr\'eductibles
remarquables de support cuspidal $({\bar L},\sigma)$ : la
repr\'esentation {\it triviale g\'en\'eralis\'ee} $\tau = \tau (\sigma_0 ,e)$ et la
repr\'esentation de {\it Steinberg g\'en\'eralis\'ee} ${\rm St}(\sigma_0
,e)$. Ce sont les seules repr\'esentations de $G$ qui apparaissent avec
multiplicit\'e  $1$ dans l'induite parabolique ${\rm Ind}_{\bar P}^{\bar
  G}\ \sigma$. La repr\'esentation de Steinberg se distingue  de la
triviale g\'en\'eralis\'ee par le fait qu'elle est g\'en\'erique.
\bigskip

 Dans [SZ], Silberger et Zink d\'eterminent les idempotents de $\HH
 ({\bar G}, {\tilde \sigma })$ correspondant \`a ces deux
 repr\'esentations. Ici c'est la triviale g\'en\'eralis\'ee qui nous
 int\'eresse.  L'idempotent associ\'ee est donn\'e par :
$$
e_\tau = \frac{1}{p_{e-1} (q_K )}\sum_{w\in W_0} {\bar f}_w\
$$

\noi o\`u $p_{e-1}$ d\'esigne le polyn\^ome de Poincar\'e du syst\`eme de racines
de type $A_{e-1}$ :
$$
p_{e-1}(X)=\prod_{k=1}^{e-1}  (1+x+\cdots + x^{k})\ .
$$

Par [SZ], page 3349, la dimension de $\tau$ est donn\'ee par
$$
{\rm dim}\ \tau ={\rm Tr}( e_\tau (1) )\vert {\bar G}\vert .
$$

\begin{proposition} Le caract\`ere de $\tau =\tau (\sigma_0
  ,e)$ en un \'el\'ement $\gamma$ de $\bar G$ est donn\'e par
$$
{\rm Tr}\ \tau (\gamma ) =\sum_{x\in {\bar G}} [{\rm Tr}\ e_\tau
](x\gamma x^{-1})\ .
$$
\end{proposition}

\noi {\it D\'emonstration}. Commen\c cons par v\'erifier que  les hypoth\`eses
2 et 3 de la section pr\'ec\'edente sont satisfaites.
\medskip

 On a $e_\tau (1) = \frac{1}{p_{e-1}(q_K)}f_1 (1)=\frac{1}{p_{e-1}(q_K
   ).\vert {\bar P}\vert }{\rm id}_{X}$, ce qui montre que $\lambda_1
 = \frac{1}{p_{e-1}(q_K ).\vert {\bar P}\vert }$ et l'hypoth\`ese 3.
\medskip

 On construit un produit scalaire $\bar P$-invariant $\langle - , - \rangle_X$
 sur $X$ de la fa\c con suivante. On fixe un produit scalaire $\langle
 -,-\rangle_0$,  ${\rm
   GL}(f,k_E )$-invariant sur $X_0$, et on pose :
$$
\langle v_1\otimes \cdots \otimes v_e , w_1 \otimes \cdots \otimes w_e
\rangle_X = \prod_{i=1,...,e} \langle v_i ,w_i \rangle_0\ .
$$
\noi En d'autres termes, on d\'efinit ce produit scalaire en d\'ecr\'etant
que si $(v_i )_i$ est une base orthonorm\'ee de $X_0$, alors
$(v_{k_1}\otimes \cdots \otimes v_{k_e})_{k_1,...,k_e}$ est une base
orthonorm\'ee de $X$.

 Alors pour $w\in W_0$,  l'adjoint
 $T_w^*$ de $T_w \in {\rm End}_\CC (X)$,  relativement \`a $\langle
 -,-\rangle_X$,  est  $T_{w^{-1}}$, et l'hypoth\`ese 2 en d\'ecoule ais\'ement.
 \medskip

 On peut donc appliquer la proposition 6.2 : pour $\gamma\in {\bar
   G}$, on a
$$
{\rm Tr}\ \tau (\gamma )
=
\frac{1}{\lambda_1}\frac{{\rm dim}\ \tau}{{\rm dim}\ \sigma}
\frac{\vert {\bar P}\vert}{\vert {\bar G} \vert}\sum_{x\in {\bar
    P}\backslash {\bar G}}{\rm Tr}(e(x\gamma
x^{-1}))\ .
$$
\noi Tenant compte du fait que
$$
{\rm dim}\ \tau ={\rm Tr}\ e_\tau (1)\ . \ \vert {\bar G}\vert =
\lambda_1 \vert {\bar G}\vert {\rm dim}\ \sigma ,
$$
\noi il vient
$$
{\rm Tr}\ \tau (\gamma )
= \vert {\bar P}\vert \sum_{x\in {\bar
    P}\backslash {\bar G}}{\rm Tr}(e(x\gamma
x^{-1})) = \sum_{x\in {\bar G}}{\rm Tr}(e(x\gamma
x^{-1}))\ .
$$

\section{Une formule de caract\`ere pour les repr\'esentations de carr\'e
  int\'egrable  de ${\rm GL}(n,F)$}

 On pose $(\pi_\lambda ,\VV_\lambda )={\mathcal E}_{\Psi} (\St_H )\in {\mathcal R}_{(J,\lambda )} (G)$, o\`u $\Psi$
est un isomorphisme unitaire d'algèbres de Hecke comme en (5.3),  et on note
 $\Theta_\lambda$ son caract\`ere d'Harish-Chandra. Notons que la classe
 d'isomorphie de
 $\pi_\lambda$,  et donc $\Theta_\lambda$, d\'epend de l'isomorphisme
 d'alg\`ebres de Hecke $\Psi$ choisi et pas seulement du type $(J,\lambda
 )$. On note $\varphi_\lambda$ le pseudo-coefficient de $\pi_\lambda$
 construit dans la section 8. Rappelons le r\'esultat fondamental
 suivant.

\begin{theorem} ([Ka], Proposition 3, page 28.) Soit $g\in G$ un
  \'el\'ement elliptique r\'egulier et soit $\mu_{G/Z}$ la mesure de Haar
  sur $G/Z$ fix\'ee dans la section 4. Alors :
$$
\Theta_\lambda (g)=\int_{G/Z} \varphi_\lambda (xg^{-1}x^{-1})\, d\mu_G
(\dot{x})\ .
$$
\end{theorem}

\begin{remark} En fait l'article de Kazhdan est \'ecrit sous les
  conditions restritives que la caract\'eristique de $F$ est nulle et
  que le centre de $G$ est compact. Cependant Badulescu (point (ii) du
  Th\'eor\`eme (4.3) de [Ba], page 64) d\'emontre ce r\'esultat pour ${\rm
    GL}(n)$ sans restriction sur le corps $F$.
\end{remark}

 L'objet de cette section est de calculer $\Theta_\lambda (g_0 )$
 lorsque $g_0$ est un \'el\'ement elliptique r\'egulier de la forme $\zeta
 u$, o\^u :
\medskip

 -- $\zeta\in U(\Bfr_M )\subset J_M$ est une racine primitive
 $q_F^N-1$ de l'unit\'e, $q_F$ \'etant la taille du corps r\'esiduel de $F$,
\smallskip

 -- $u\in H^{1}(\beta ,\Afr_M )\subset J_M$, o\^u $H^1 (\beta ,\Afr_M )$
 est le groupe d\'efini dans [BK](3.1).

\begin{lemma} (i) L'image $\bar{g}_0 =\bar{\zeta}$ de $g_0$ dans $J_M
  /J_M^1 \simeq U(\Bfr_M)/U^1 (\Bfr_M )\simeq {\rm GL}(ef,k_E )$ est
  elliptique r\'eguli\`ere.
\smallskip

\noi (ii) Plus g\'en\'eralement, si $g\in U(\Afr_M )$ est tel que
$g^{-1}g_0 g\in J_M$,  alors l'image de $g^{-1}g_0 g$ dans $J_M /J_M^1
\simeq {\rm GL}(ef,k_E )$ est elliptique r\'eguli\`ere.
\smallskip

\noi (iii) L'\'el\'ement $g_0$ est minimal sur $F$ au sens de
     [BK](1.4.14), page 41, et $F[g_0 ]/F$ est une extension non
     ramifi\'ee de degr\'e maximal $N$. L'ordre $\Afr_M$ est l'unique
     ordre h\'er\'editaire  de $A$ normalis\'e par $g_0$. 
\end{lemma}

\noi {\it D\'emonstration}. Le point (i) est facile. Montrons d'abord
(iii). Il d\'ecoule imm\'ediatement de la d\'efinition de [BK] que $\zeta$
est minimal. La strate $[\Afr_M, 0,-1 ,\zeta ]$ est donc
simple (voir [BK](1.5) pour plus de détails). Puisque 
$$
\zeta u -\zeta =\zeta (u-1)\in {\rm Rad}\, (\Bfr_M ) \subset {\rm
  Rad}\, (\Afr_M )\ ,
$$
\noi on a l'\'equivalence de strates :
$$
[\Afr_M ,0,-1 ,\zeta u]\sim [\Afr_M ,0,-1 ,\zeta ]\ .
$$
\noi Donc par la proposition (2.2.2) de [BK], page 52, la strate
$[\Afr_M ,0,-1,\zeta u]$ est simple et l'extension $F[\zeta u]/F$ est
de degr\'e maximal. En particulier $\zeta u$ est minimal sur $F$. Par
[BK](2.1.4), page 50,  $F[\zeta u]/F$ est non ramifi\'ee. Enfin il est
clair que $\Afr_M$ est normalis\'e par $g_0$ et par l'exercice (1.5.6)
de [BK], page 44, c'est l'unique ordre h\'er\'editaire de $A$ ayant cette
propri\'et\'e.

 Pour prouver (ii), supposons par l'absurde que l'image de
 $gg_0g^{-1}$ ne soit pas elliptique r\'eguli\`ere. Il existerait alors un
 $\ofr_E$-ordre h\'er\'editaire $\Bfr_M '$ strictement inclus dans
 $\Bfr_M$ tel que $gg_0 g^{-1}\in U(\Bfr_M ' )J_M^1$. Ainsi on aurait
$$
gg_0 g^{-1}\in U(\Afr_M ')J_M^1 \subset U(\Afr_M ')U^1 (\Afr_M
)\subset U(\Afr_M ')
$$
\noi o\^u $\Afr_M '$ est l'$\ofr_F$-ordre h\'er\'editaire de $A$
correspondant \`a $\Bfr_M '$ par la correspondence de [BK](1.2). On
obtiendrait alors $g_0 \subset U(g^{-1}\Afr_M 'g)$, avec
$g^{-1}\Afr_M g$ sous-ordre h\'er\'editaire strict de $\Afr$, ce qui
contredirait la minimalit\'e de $g_0$.
\bigskip

 Soit $g\in G$, nous avons vu dans la proposition (8.3) que $\varphi
 (gg_0^{-1}g^{-1})$ est nul \`a moins que $gg_{0}^{-1}g^{-1}\in
 J_M$. D'un autre c\^ot\'e, si $gg_0^{-1}g^{-1}\in J_M \subset U(\Afr_M
 )$, alors, par minimalit\'e de $g_0$, $g$ appartient au normalisateur
 de $\Afr_M$, c'est-\`a-dire $F^\times U(\Afr_M )$. Pour d\'eterminer la
 fonction $g\mapsto \varphi_\lambda (gg_0^{-1}g^{-1})$, nous supposons
 donc que $g\in F^\times U(\Afr_M )$ et $gg_0^{-1}g^{-1}\in J_M$. Par
 le point (ii) du lemme pr\'ec\'edent, on a que l'image de $gg_0^{-1}g^{-1}$
 dans $J_M /J_M^1$ est elliptique r\'eguli\`ere.

 Par la proposition (8.3), nous avons :
\begin{equation}
\varphi (gg_0^{-1}g^{-1})=\frac{(-1)^{e-1}}{e'}\sum_{T\subset
  S}\frac{(-1)^{d_T}}{(d_S +1)\mu_H (P_T )}\sum_{w\in T} {\rm Tr}\,
\tilde{\kappa}_M (gg_0^{-1}g^{-1}).{\rm Tr}\, \bar{f}_w
(\overline{gg_0^{-1}g^{-1}})\ .
\end{equation}
 
 Soient $T\subset S$ et $w\in \langle T\rangle$. Alors si $\bar{f}_w
 (\overline{gg_0^{-1}g^{-1}})\not= 0$, l'\'el\'ement elliptique
 r\'egulier $\overline{gg_0^{-1}g^{-1}}$ appartient \`a
 $\bar{P}\langle T\rangle \bar{P}$. Or si $T$ est strictement contenu
 dans $S$, ce sous-ensemble est contenu dans un parabolique propre de
 $\bar{G}$, ce qui n'est pas possible. On en d\'eduit l'expression
 suivante :

\begin{equation}
\varphi (gg_0^{-1}g^{-1})=\frac{(-1)^{e-1}}{e'} \frac{(-1)^{d_S}}{(d_S
  +1)\mu_H (P_S )}\sum_{w\in W_0}
{\rm Tr}\,\tilde{\kappa}_M (gg_0^{-1}g^{-1}).{\rm Tr}\, \bar{f}_w
(\overline{gg_0^{-1}g^{-1}})\ .
\end{equation}

Soit $\tau = \tau (\check{\sigma}_0,e )$ la repr\'esentation de
triviale g\'en\'eralis\'ee de $\bar G$ attach\'ee \`a $\check{\sigma}_0$. Alors
d'apr\`es les r\'esultats rappel\'es en section 10, elle admet pour
idempotent
$$
e_\tau = \frac{1}{p_{e-1}(q_K )} \sum_{w\in W_0} \bar{f}_w\ .
$$

Nous obtenons donc :

\begin{equation}
\varphi (gg_0^{-1}g^{-1})=\frac{(-1)^{e-1}}{e'} \frac{(-1)^{d_S}}{(d_S
  +1)\mu_H (P_S )}p_{e-1}(q_K )
{\rm Tr}\,\tilde{\kappa}_M (gg_0^{-1}g^{-1}).{\rm Tr}\, e_\tau
(\overline{gg_0^{-1}g^{-1}})\ .
\end{equation}
  
 Pour simplifier les notations, posons :
\begin{equation}
 C_S = \frac{(-1)^{e-1}}{e'} \frac{(-1)^{d_S}}{(d_S
  +1)\mu_H (P_S )}p_{e-1}(q_K )\ .
\end{equation}

En utilisant la proposition (10.1)  reliant la trace de la repr\'esentation $\tau$ \`a
l'idempotent $e_\tau$, on obtient

$$
\int_{J_M} \varphi_\lambda (gg_0^{-1}g^{-1})\, d\mu_G (g) 
$$

\begin{eqnarray}
=& C_S  \mu_G (J_M^1 ) {\rm Tr}\, \tilde{\kappa}_M (gg_0^{-1}g^{-1})\, {\rm
  Tr}\, \tau (\overline{gg_0^{-1}g^{-1}})\\
 =&  C_S \mu_G (J_M^1 ) {\rm Tr}\, \kappa_M (g^{-1}g_0g)\, {\rm
  Tr}\, \tau (\sigma_0,e) (\overline{g^{-1}g_0^{-1}g})        \\
 =&   C_S \mu_G (J_M^1 ) {\rm Tr}\, \lbrace \kappa_M \otimes\tau
 (\sigma_0 ,e)\rbrace (g^{-1}g_0^{-1}g)  \ .
\end{eqnarray}

\noi o\^u on a utilis\'e le fait que la contragr\'ediente de $\tau$ est
$\tau (\sigma_0 ,e)$.

\begin{proposition} (a)  Les repr\'esentations $\tau (\sigma_0 ,e)$
 et ${\rm St}(\sigma_0 ,e)$ sont images l'une de l'autre 
dans la dualit\'e d'Alvis-Curtis.

\noi (b)  Pour tout \'el\'ement elliptique r\'egulier 
$\alpha$  de $\bar G$, on a :
$$
{\rm Tr}\  \tau (\sigma_0 ,e)(\alpha )=(-1)^{e-1}
 {\rm Tr}\  {\rm St}(\sigma_0 ,e)(\alpha)\ .
$$
\end{proposition}

\noi {\it D\'emonstration}. Pour le point a), nous renvoyons
 \`a [DM] Corollaire 14.47.
\medskip

 La dualit\'e d'Alvis-Curtis, originellement d\'efinie au niveau 
des caract\`eres des repr\'esentations, a \'et\'e explicitement
construite par Deligne et Lusztig [DL].  Leur Corollaire (c)
 du paragraphe 5, page 290 affirme que si $E$ est une repr\'esentation
irr\'eductible de $\bar G$ de duale $E^\sharp$, alors on a 
l'\'egalit\'e suivante dans le groupe de Grothendieck des $\bar G$-modules
virtuels :
$$
(-1)^{i_0}\ E^\sharp = \sum_{I\subset \bar S} (-1)^{\vert I\vert } E_{(I)}\ ,
$$
\noi o\`u :
\bigskip

 -- $\bar S$ est le syst\`eme g\'en\'erateur d'involutions de la 
$BN$-paire sph\'erique de $\bar G$,
\medskip

 -- $E_{({\bar S})}=E$,
\medskip

 -- si $I$ est un sous-ensemble strict de $\bar S$, $E_{(I)}$ est une 
repr\'esentation de $\bar G$ induite \`a partir d'un
parabolique {\it propre},
\medskip

 -- $i_0 ={\rm Min}\ \{ \vert I\vert\ ; \ E^{U_I}\not= 0\}$, o\`u
 $U_I$ 
est le radical unipotent du parabolique standard de type $I$.
\bigskip

 Appliquons ceci \`a $E={\rm St}(\sigma_0 ,e)$. Puisque le support
 cuspidal de cette repr\'esentation est $({\bar L}, \sigma)$, on a
$i_0 = e(f-1)$. De plus si $\alpha$ est elliptique r\'egulier et $I$ 
un sous-ensemble strict de $\bar S$, on a ${\rm Tr}\ E_{(I)}
(\alpha )=0$. On obtient ainsi
$$
(-1)^{e(f-1)}{\rm Tr}\ E^{\sharp}(\alpha )=(-1)^{\vert {\bar S}\vert} 
{\rm Tr}\ {\rm St}(\sigma_0 ,e)(\alpha)=(-1)^{ef-1}{\rm Tr}\
{\rm St}(\sigma_0 ,e)(\alpha)
$$
\noi et l'assertion b) en d\'ecoule.
\bigskip

 Pour r\'esumer, nous avons montr\'e le r\'esultat suivant.
\bigskip

\begin{proposition} Avec les notations pr\'ec\'edentes, on a :
\begin{equation}
\int_{J_M}\varphi_\lambda (gg_0^{-1}g^{-1}) d\mu_G (g)=
C_S \mu_G (J_M^1 ) (-1)^{e-1}
\ {\rm Tr}\ \lbrace \kappa_M \otimes {\rm St}(\sigma_0 ,e)\rbrace ( g_0 )\ .
\end{equation}
\end{proposition}

\bigskip

 Ecrivons 

$$
\int_{U(\Afr_M )} \varphi_\lambda (gg_0^{-1}g^{-1})\, d\mu_G
(g)  =  \sum_{u\in J_M\backslash U(\Afr_M )}
\int_{J_M} \varphi_\lambda  (gu\gamma u^{-1}g^{-1})\, d\mu_G (x)
$$
$$ =  \sum_{u\in J_M\backslash U(\Afr_M )} C_S \mu_G (J_M^1
) (-1)^{e-1} {\rm Tr}\, \lbrace \kappa_M\otimes {\rm St}(\sigma_0 ,e)\rbrace (ug_0 u^{-1}) 
$$
\noi c'est-\`a dire :
\begin{equation}
\int_{U(\Afr_M )} \varphi_\lambda (gg_0^{-1}g^{-1})\, d\mu_G
(g) =  \sum_{u\in U(\Afr_M )/J_M } C_S \mu_G (J_M^1
) (-1)^{e-1} {\rm Tr}\, \lbrace \kappa_M\otimes {\rm St}(\sigma_0 ,e)\rbrace (u^{-1}g_0
u) 
\ .
\end{equation}

La fonction 
$$
x\mapsto \varphi_\lambda (x\gamma^{-1}x^{-1})
$$
\noi \'etant \`a support dans $F^{\times} U(\Afr_M )$ et \'etant invariante
 par $F^\times$,  le membre de gauche de notre derni\`ere \'equation
 est en r\'ealit\'e \'egal \`a 
$$
\int_{G/Z} \varphi_\lambda (gg_0^{-1}g^{-1})\, d\mu_{G/Z}(\dot{g})\ .
$$

Nous allons simplifier le membre de droite. On a 
$$
C_S \mu_G (J_M^1 ) (-1)^{e-1} =\frac{[(-1)^{e-1}]^2}{e'}
\frac{(-1)^{d_S}}{(d_S +1)\mu_H (P_S )}
p_{e-1}(q_K )\mu_G (J_M^1 )\ ,
$$
\noi avec
\medskip

 -- $d_S =0$ (le simplexe correspondant \`a $S$ est un sommet),
\smallskip

 -- $\mu_G (J_M^1 )=1$ (normalisation de la mesure de Haar sur $G$),
\smallskip

 -- $e' =e(E/F)=1$, car $E/F$ est non ramifi\'ee,
\smallskip

 -- $P_S ={\rm GL}(e,\ofr_K )$.
\medskip

De plus, puisque $\mu_H (I)=1$, on a 
$$
\mu_H (P_S )=\vert P_S /I\vert = \vert {\rm GL}(e,k_K )/B_e (k_K )\vert
$$
\noi o\^u $B_e$ est le sous-groupe de Borel standard de ${\rm GL}(e)$. 
Un calcul classique donne donc
$\mu_H (P_S )= p_{e-1}(q_K )$.

\begin{remark} Il est amusant de noter que $p_{e-1}(q_K )$ est la 
$q_K$-factorielle de $e$, qui se trouve \^etre le nombre de points
rationnels de la vari\'et\'es de drapeaux complets sur $k_K^e$ ...
\end{remark}

 Nous venons donc de d\'emontrer la formule.  de caract\`ere suivante

\begin{theorem} Avec les notations pr\'ec\'edentes, la valeur du caract\`ere
 de Harish-Chandra 
de $\pi_\lambda$ en $\zeta u $ est donn\'ee par
\begin{equation}
\Theta_\lambda (\zeta u)=\sum_{v\in U(\Afr_M)/J_M} {\rm Tr}\,
 \lbrace \kappa_M\otimes {\rm St}(\sigma_0 ,e)\rbrace (v^{-1}(\zeta u)
v) 
\ .
\end{equation}
\end{theorem}

\section{Transfert du pseudo-coefficient de Kottwitz-II}

 Avec les notations de la section 1, 
nous faisons ici les deux hypoth\`eses suivantes :
\medskip

(i) l'extension $E/F$ est totalement ramifi\'ee,
\smallskip

 (ii) l'ordre $\Bfr$ est minimal.
\medskip

 Il s'ensuit que la repr\'esentation $\sigma$ intervenant dans le type
 $\lambda =\kappa\otimes \sigma$ est un caract\`ere de $U(\Bfr )/U^1
 (\Bfr )\simeq (k_E^{\times})^{N/[E:F]}$. Quitte \`a modifier la
 b\^eta-extension $\kappa$, nous ne perdons donc rien \`a supposer que
 $\sigma$ est le caract\`ere trivial et $\lambda =\kappa$.
\medskip

 Notons qu'ici on a  $K=E$ et le groupe $H$ est ${\rm
   GL}(N/[E:F],E)$. L'ordre $\Afr$ est minimal. Fixons une
 uniformisante $\varpi_E$ de $E$ et une 
 uniformisante $\Pi =\Pi_\Bfr$ de l'ordre $\Bfr$, choisie de telle
 sorte que $\Pi^{N/[E:F]} =\varpi_E$ ; $\Pi$ est aussi une
 uniformisante de l'ordre $\Afr$.
\bigskip 

 Fixons un isomorphisme unitaire d'alg\`ebres de Hecke :
$$
\Psi~: \ \HH (H,{\mathbf 1}_I )\lra \HH (G,\lambda )
$$
\noi et notons
$$
{\mathcal E}_\Psi~: \ \Rep_{(I,{\mathbf 1}_I)} \lra \Rep_\lambda (G)
$$
\noi l'\'equivalence de cat\'egories correspondante. Finalement posons
$$
\pi_\lambda ={\mathcal E}_\Psi (\St_H )\ 
$$
\noi et notons $\omega$ le caract\`ere central de $\pi_\lambda$.
\bigskip

 Transf\'erons le pseudo-coefficient de Kottwitz $f_0\in \HH (H,{\mathbf
   1}_I )$ en un \'el\'ement $\varphi_\lambda$ de $\HH (G, \lambda\omega
 )$ par les \'etapes suivantes :
\medskip

 (i) $F_\lambda =\Psi (F_0 )\in \HH (G,\lambda )$,
\smallskip

 (ii) $f_\lambda =P_{\omega} (F_\lambda )\in \HH (G,\lambda \omega )$,
\smallskip

 (iii) $\varphi_\lambda ={\rm Tr}_{W_{\check{\lambda}}}\circ f_\lambda
\in e_{\lambda\omega}\star \HH (G)\star e_{\lambda \omega}$.
\medskip

 Nous avons choisi ici la mesure de Haar $\mu_G$  sur $G$ telle que
 $\mu_G (J)= 1$, la mesure de Haar sur $F^\times$ normalis\'ee par
 $\mu_{F^\times} (\ofr^\times )=1$, et la mesure de Haar quotient
 $\mu_{G/Z}$ sur $G/Z$. La projection $P_\omega$ et l'idempotent
   $e_{\lambda\omega}$ sont alors d\'efinis de fa\c con habituelle.
\bigskip

 Le r\'esultat suivant se d\'emontre comme la proposition (6.5).

\begin{proposition} La fonction $\varphi_\lambda$ est un
  pseudo-coefficient de $\pi_\lambda$.
\end{proposition}

 Notons $\tilde{E} = E[\Pi ]=F[\Pi ]$. C'est une extension totalement
 ramifi\'ee de $F$ de degr\'e $N$ dont le groupe multiplicatif normalise
 $\Bfr$ et $\Afr$. En fait $\Pi$ est minimal sur $F$ au sens de
 [BK](1.4.14) et $\Afr$ est l'unique ordre h\'er\'editaire de $A$
 normalis\'e par $\Pi$ (cf. [BK], exercise (1.5.6)).
\bigskip

 Nous allons d\'eterminer les valeurs $\varphi (xg_0 g^{-1})$, $x\in G$,
 pour un \'el\'ement $g_0$ de la forme :
$$
g_0 =\Pi^{-\nu} u , \ \nu\in \{ 0,1,...,N-1\} ,\ {\rm pgcd}(\nu
,N)=1\ \text{et} \ u\in H^{1}(\beta
,\Afr )\ 
$$

\begin{lemma} L'\'el\'ement $g_0$ est minimal sur $F$ et normalise
  $\Afr$. En particulier $\Afr$ est l'unique ordre h\'er\'editaire de $A$
  normalis\'e par $g_0$. 
\end{lemma}

\noi {\it D\'emonstration}. On a $u\in 1+{\rm Rad}(\Afr )\subset U(\Afr
)$ et $g_0$ normalise donc $\Afr$.  Les strates $[\Afr ,-\nu ,-\nu
  -1,g_0 ]$ et $[\Afr ,-\nu ,-\nu -1 ,\Pi^{-\nu}]$ sont \'equivalentes
et on conclut comme dans le lemme (11.4)(ii).
\bigskip

 On a, rappelons-le :
$$
\varphi_\lambda (g)=
$$
$$
\frac{(-1)^{e-1}}{e'}\sum_{k\in \ZZ}\sum_{T\subset S}\frac{(-1)^{d_T}}{(d_T
  +1)\mu_H (P_T )} \sum_{w\in \langle T\rangle} \sum_{l=0}^{e'n_T
  -1}\epsilon_T^{l}\omega (\varpi_F )^{-k} {\rm
  Tr}_{W_{\check{\lambda}}}\circ \Psi (f_{z_T^l w})(\varpi_F^{-k}g)\,
$$
\noi pour $g\in G$ et o\^u $e'=e(E/F)$.
\medskip

 Si $g=xg_0^{-1}x^{-1}$, alors $v_F ({\rm det}\, (g))=\nu$. D'un autre
 c\^ot\'e la fonction $\Psi (f_{z_T^l w})$ a un support contenu dans 
$J z_T^l w J$. Or tout \'el\'ement $y$ de $Jz_T^l wJ$ v\'erifie
$$
v_F ({\rm det}\, (y))=lv_F ({\rm det}\, (z_T ))=l\frac{N}{n_T e'}\ .
$$
\noi Ainsi,  si le terme $\Psi (f_{z_T^l w})(\varpi^{-k}g)$ est non
nul, on a $\ds l\frac{N}{n_T e'} =\nu -kN$. Puisque $0\leq
l\frac{N}{n_T e'} <N$, ceci entra\^ine $k=0$. Donc $ l\frac{N}{n_T
  e'}=\nu$, c'est-\`a-dire $\nu =l$, puisque $\nu$ est premier \`a
$N$. Mais alors $n_T =N/e' =N/[E:F]$, ce qui entra\^ine $T=\emptyset$. 

 Pour r\'esumer, si $xg_0^{-1}x^{-1}$ est dans le support de
 $\varphi_\lambda$, on a 
$$
\varphi_\lambda (xg_0^{-1}x^{-1})=\frac{(-1)^{e-1}}{e'}
\frac{(-1)^{d_\emptyset}}{(d_\emptyset +1)\mu_H
    (I)}\epsilon_{\emptyset}^\nu {\rm Tr}_{W_{\check{\lambda}}}\circ
  \Psi (f_{\Pi^\nu})(xg_{0}^{-1}x^{-1})\ .
$$
Notons aussi que $d_\emptyset =N/e' -1$ et $e=N/e'$. De plus $g_0$
\'etant minimal sur $F$ et $\Psi (f_{\Pi^\nu})$ de support contenu dans
${\rm N}_G (\Afr ) =\langle \Pi \rangle\, U(\Afr )$, on a $
\varphi_\lambda (xg_0^{-1}x^{-1})=0$, si $x\not\in {\rm N}_G (\Afr )$.
Nous avons donc d\'emontr\'e le

\begin{lemma} On a 
$$
\varphi_\lambda (xg_{0}^{-1}x^{-1}) =\frac{\epsilon_\emptyset}{N}{\rm
  Tr}_{W_{\check{\lambda}}}\circ \Psi (f_\Pi )^\nu (xg_{0}^{-1}x^{-1}) , \ x\in N_G (\Afr )
$$
\noi et 
$$
\varphi_\lambda (xg_{0}^{-1}x^{-1})=0\text{ si } x\not\in {\rm
  N}_{G}(\Afr )\ .
$$
\end{lemma}

 Notons que $\epsilon_\emptyset$ est la signature de la permutation
 circulaire
$$
\left(
\begin{array}{lllll}
1 & 2 & \cdots & n-1 & n \\
n & 1 & \cdots & n-2 & n-1
\end{array}
\right)
$$
\noi c'est-\`a-dire $(-1)^{n-1}$.

\bigskip

 Notons $V$ l'espace de $\pi_\lambda$

\begin{proposition} (i) La composante $\lambda$-isotypique $V^\lambda$
  est isomorphe \`a $\lambda$ comme repr\'esentation de $J$
  (i.e. $\lambda$ intervient dans $\pi$ avec multiplicit\'e $1$).
\smallskip

\noi (ii) L'espace  $V^\lambda$ est stable sous
  l'action de $N_{B^\times}(\Bfr )$. En particulier l'action de
  $\langle \Pi \rangle J$ dans $V^\lambda$ est une repr\'esentation
  $\tilde \lambda$ qui prolonge la repr\'esentation $\lambda$.
\smallskip

\noi (iii) Il existe une constante non nulle $c\in \CC$ telle que
$$
\Psi (f_\Pi )(\Pi j)=c.{\tilde \lambda}(\Pi j)\ , \ j\in J\ .
$$

\noi (iv) Pour $k\geqslant 0$, on a 
$$
\Psi (f_\Pi )^k (\Pi^k j)=c^k .{\tilde \lambda}(\Pi^k  j)\ , \ j\in J\ .
$$
\end{proposition}

\noi {\it D\'emonstration}. Le $\HH (H,I)$-module correspondant \`a $\St_H$
est de dimension $1$. Il en est donc de m\^eme du $\HH (G,\lambda
)$-module $M$ correspondant \`a $V$. Comme $\CC$-espace vectoriel la
composante isotypique $V^\lambda$ est isomorphe \`a
$M\otimes_{\CC}W_{\check{\lambda}}$. On a donc ${\rm dim}_\CC\,
(V^\lambda )={\rm dim}\, (\lambda )$, ce qui prouve (i).

Par la
Proposition (5.5.11), page 185 de [BK], l'entrelacement de $(J,\lambda
)$ est $JB^\times J$, en particulier il contient $N_{B^\times}(\Bfr
)$. Mais par (3.1.15)(ii), le groupe $N_{B^\times}(\Bfr )$ normalise
$J$ et il normalise donc la paire $(J,\lambda )$. Le point (ii) en d\'ecoule. 

 Par [BK](5.66)(i), page 190, toutes les fonctions de $\HH (G,\lambda
 )$ \`a support dans $\Pi .J$ sont proportionnelles. Par le point (i),
 on a que $\Pi$ entrelace $\lambda$ avec comme op\'erateur
 d'entrelacement ${\tilde \lambda}(\Pi )$. Donc il existe une fonction
 $f\in \HH (G,\lambda )$, non nulle \`a support $\Pi .J$, telle que
 $f(\Pi )={\tilde \lambda}(\Pi )$. On en d\'eduit qu'il existe $c\in
 \CC^\times$ tel que $\Pi (f_{\Pi})$, qui est \`a support $\Pi .J$, vaut
 $c.{\tilde \lambda}$ en $\Pi$. Le point (iii) en d\'ecoule.  Le (iv) se
 d\'eduit de (iii), par r\'ecurrence,  via un calcul imm\'ediat de convolution (utilisant le
 fait que $\mu_G (J)=1$).

 \begin{theorem} Avec les notations pr\'ec\'edentes, supposons que le
polyn\^ome caract\'eristique de $g_0 = \Pi^{-\nu}u\in A$ est s\'eparable sur
$F$, de sorte que $g_0$ est elliptique r\'egulier. Le caract\`ere d'Harish-Chandra
$\Theta_\lambda$ de $\pi_\lambda$ en $g_0$ est alors donn\'e par la
formule
$$
\Theta_\lambda (\Pi^{-\nu}u)= (-1)^{\nu (n-1)}\, c^\nu \sum_{x\in
  U(\Afr )/J} {\rm Tr}\, {\tilde \lambda}\, \big( x^{-1}(\Pi^{-\nu
  }u)x\big)\ .
$$
\end{theorem}

\noi {\it D\'emonstration}. En \'etendant la fonction ${\rm Tr}\, {\tilde
  \lambda}$ par $0$ \`a $G$ tout entier, on a 
$$
\varphi_\lambda  (xg_0^{-1}x^{-1}) =\frac{(-1)^{\nu (n-1)}}{N} {\rm
  Tr}\, {\tilde \lambda} (x^{-1}(\Pi^{-\nu}u)x)\ , \ x\in G\ .
$$

\noi La formule de Kazhdan (Th\'eor\`eme (11.2)) donne alors
successivement :

\begin{eqnarray}
\Theta_\lambda (g_0 ) & = & \frac{(-1)^{\nu (n-1)}}{N}
\sum_{\langle \Pi\rangle U(\Afr )/F^\times} {\rm Tr}\, {\tilde
  \lambda}(x^{-1}g_0 x)\, d\mu_{G/Z}(\dot{x})\\
        &=&  \frac{(-1)^{\nu (n-1)}}{N}\sum_{k=0}^{N}\int_{\Pi^k
  F^\times U(\Afr )/F^\times} {\rm Tr}\, {\tilde
  \lambda}(x^{-1}g_0 x)\, d\mu_{G/Z}(\dot{x})\\
  & = & \frac{1}{N}\sum_{k=0}^{N-1} (-1)^{\nu (n-1)}c^\nu \int_{\Pi^k
  U(\Afr )}  {\rm Tr}\, {\tilde \lambda}(x^{-1}g_0 x)\, d\mu_G (x)\\
 &=& c^{\nu}(-1)^{\nu (n-1)} \int_{U(\Afr )}{\rm Tr}\, {\tilde
  \lambda}(x^{-1}g_0 x)\, d\mu_G (x)\\
 &=& c^{\nu}(-1)^{\nu (n-1)} \sum_{x\in U(\Afr )/J} {\rm Tr}\,
  {\tilde \lambda}(x^{-1}g_0  x)\ .
\end{eqnarray}

\appendix

\section{Quelques lemmes sur l'action d'une alg\`ebre de Hecke
  sph\'erique}

 Les r\'esultat de cette annexe sont peut \^etre bien connus, mais
 l'auteur, ne connaissant pas de r\'ef\'erences, a pr\'ef\'er\'e les d\'emontrer.
\bigskip
 
 On fixe~:
\medskip

 -- un groupe $G$ localement profini et unimodulaire et une
 mesure de Haar $\mu$ sur $G$ ; 
\smallskip

 -- une paire $(J,\lambda )$ form\'ee d'un sous-groupe ouvert compact
 $J$ de $G$ et d'une repr\'esentation lisse irr\'eductible  $(\lambda ,
 W)$ de $J$.
\medskip

 Soit $(\check{\lambda}, \check{W})$ la repr\'esentation contragr\'ediente
 de $(\lambda , W )$. Rappelons que l'alg\`ebre de Hecke sph\'erique $\HH
 (G,\check{\lambda})$ est le $\CC$-espace vectoriel form\'e des fonctions
 $$
\varphi ~: \ G\lra {\rm End}_\CC\ (W)
$$
\noi qui sont \`a support compact et se transforment selon
$$
\varphi (j_1 gj_2 )=\lambda (j_1 )\circ  \varphi (g)\circ \lambda (j_2 )\ ,
$$
\noi muni du produit de convolution
$$
\varphi_1 \star \varphi_2 (g)=\int_G \varphi_1 (x)\circ 
\varphi_2 (x^{-1}g)\ dx\ , g\in G\ .
$$

 Il est classique qu'on a un isomorphisme canonique d'alg\`ebres
$$
\HH (G,\check{\lambda})\lra {\rm End}_G\, (\cind{J}{G}\ \lambda )
$$
\noi qui \`a $\varphi\in \HH (G,\check{\lambda})$ associe l'op\'erateur de
convolution
$$
f\mapsto \varphi\star f\ , \ f\in \cind{J}{G}\ \lambda\ ,
$$
\noi o\^u
$$
\varphi\star f(g)=\int_G \varphi (x)[f(x^{-1}g)]\ dx\ g\in G\ .
$$
\noi On notera $\varphi_\star$ l'op\'erateur de convolution attach\'e \`a
$\varphi$. 
\medskip

 Fixons \`a pr\'esent une repr\'esentation lisse $(\pi ,\VV )$ de $G$. On
 lui associe deux espaces vectoriels :
$$
M_\VV ={\rm Hom}_J\ (\lambda ,\pi )\text{ et } {\tilde M}_\VV
={\rm Hom}_G \ (\cind{J}{G}\ \lambda ,\pi )\ .
$$
\noi Par r\'eciprocit\'e de Frobenius pour l'induction compacte, ces deux
espaces sont canoniquement isomorphes. Rappelons comment est r\'ealis\'e
cet isomorphisme canonique. Consid\'erons l'\'el\'ement $\ds \alpha\in {\rm
  Hom}_J (\lambda , \cind{J}{G}\ \lambda )$ qui \`a $w\in W$ associe la
fonction $T_w$ \`a support dans $J$ donn\'ee par $T_w (j)=\lambda (j).w$,
$j\in J$. L'isomorphisme canonique $\Psi$~: ${\tilde M}_\VV \lra
M_\VV$ est donn\'e par $\Psi ({\tilde \varphi})={\tilde \varphi}\circ
\alpha$, ${\tilde \varphi}\in {\tilde M}_\VV$. Nous avons besoin
d'exhiber son inverse.

\begin{lemma} L'inverse $\Phi$~: $M_\VV \lra {\tilde M}_\VV$ de $\Psi$
  est donn\'e par $\Phi (\varphi )={\tilde \varphi}$, o\^u 
$$
{\tilde \varphi}(f) =\frac{1}{\mu (J)}\ \int_G \pi(x).\varphi
(f(x^{-1}))\ dx\ , \ f\in \cind{J}{G}\ \lambda\ .
$$
\end{lemma}

\noi {\it D\'emonstration}. Soient $\varphi\in M_\VV$ et $w\in W$. On a
successivement :

\begin{eqnarray}
\Psi (\Phi (\varphi ))(w) & = & \Phi (\varphi )\circ \alpha\ (w) \\ 
      & = &  \Phi (\varphi )(T_w ) \\
      & = & \frac{1}{\mu (J)}\ \int_G \pi (x).\varphi (T_w (x^{-1}))\,
dx\\
      & = & \frac{1}{\mu (J)}\ \int_{J} \pi (j). \varphi (\lambda
(j^{-1}).w)\, dj\\
      &= & \frac{1}{\mu (J)}\ \int_{J} \pi (j)\circ \pi (j^{-1}).\varphi
(w)\, dj\\
      & = &  \varphi (w)\ .
\end{eqnarray}

\noi On a bien $\Psi (\Phi (\varphi )) =\varphi$.
\bigskip

On regarde ${\tilde M}_{\VV}$ comme un $\HH
(G,\check{\lambda})$-module \`a droite via
$$
{\tilde \varphi}.f= {\tilde \varphi}\circ f_\star\ , \ {\tilde
  \varphi}\in {\tilde M}_\VV\ , \ f\in \HH (G,\check{\lambda})\ .
$$
\noi On munit de m\^eme $M_\VV$ d'une structure de $\HH
(G,\check{\lambda})$-module \`a droite {\it via} l'identification
canonique $M_\VV \simeq {\tilde M}_\VV$, c'est-\`a-dire :
$$
\varphi .f=\Psi (\, \Phi (\varphi )\circ f_\star \, )\ , \ \varphi\in
M_\VV\ , \ f\in  \HH (G,\check{\lambda})\ .
$$
\noi Nous allons calculer cette action de fa\c con explicite.

\begin{lemma} L'action naturelle de $\HH (G,\check{\lambda})$ sur
  $M_\lambda$ est donn\'ee par
$$
\varphi .f = \int_G \pi (x)\circ \varphi\circ f(x^{-1})\, dx\ , \ \varphi\in
M_\VV\ , \ f\in  \HH (G,\check{\lambda})\ .
$$
\end{lemma}

\noi {\it D\'emonstration}. Soient $\varphi\in {\rm Hom}_J\ (\lambda
,\pi )$, $f\in \HH (G,\check{\lambda})$ et $w\in W$. On a
successivement :

\begin{eqnarray}
(\varphi .f )(w) & =& \Psi (\, \Phi (\varphi )\circ f_\star\, )(w)\\
 & = & (\, \Phi (\varphi\circ f_\star \, ) (T_w)\\
 & = & \Phi (\varphi (\, f\star T_w \, )\\
 & = & \frac{1}{\mu (J)}\int_G \pi (x) . \varphi \{ (f\star T_w
  )(x^{-1})\}\, dx\\
 & = &   \frac{1}{\mu (J)}\int_G \pi (x) . \varphi \big\{ \, \int_G
  f(u)[\, T_w (u^{-1}x^{-1})\, ]\, du\, \big\}\, dx\\
 & = &    \frac{1}{\mu (J)}\int_G \pi (x) . \varphi \big\{\, \int_J
  f(x^{-1}j^{-1})[\, T_w (j)\, ]\, dj\big\}\, dx\\
 & = &   \frac{1}{\mu (J)}\int_G \pi (x) . \varphi \big\{\, \int_J
f(x^{-1})\circ \lambda (j^{-1})\circ \lambda (j)(w)\, dj\big\}\, dx\\
 & = &  \frac{1}{\mu (J)}\int_G \pi (x) . \varphi \big\{\,  \mu
(J).f(x^{-1})(w)\, \big\}\, dx\\
 & = & \big( \int_G \pi (x)\circ \varphi \circ f(x^{-1})\big) (w)
\end{eqnarray}

\noi \noi o\^u \`a la ligne (6), on a fait le changement de variable
$j=u^{-1}x^{-1}$. Le r\'esultat en d\'ecoule.

\bigskip

 Nous allons \'enoncer des r\'esultats similaires en rel\^achant les
 hypoth\`eses. Les d\'emonstrations sont laiss\'ees au lecteur. Nous gardons
 donc les m\^emes notations mais modifions les hypoth\`eses de la fa\c con
 suivante :
\medskip

-- $J$ est ouvert, contient le centre $Z$ de $G$ et $J/Z$ est compact
; 

\smallskip

 -- $(\pi ,\VV )$ poss\`ede un caract\`ere central $\omega_\pi$ qui
 co\"incide avec le caract\`ere central de $\lambda$.
\medskip

 On fixe une mesure de Haar $\mu_Z$ sur $Z$ et on note $\mu_{G/Z}$ la
 mesure quotient de $\mu$ par $\mu_Z$.  Cette mesure de Haar d\'efinit
 une alg\`ebre de Hecke $\HH (G,\check{\lambda} )$. Cette derni\`ere s'identifie \`a
 l'alg\`ebre ${\rm End}_G \, (\cind{J}{G}\, \lambda )$. Cette
 identificationn envoie $\varphi\in \HH (G,\check{\lambda})$ sur l'op\'erateur
 $\varphi_\star$ donn\'e par
$$
\varphi_\star (f)(g)=f\star \varphi (g)=\int_{G/Z}\varphi (x)
       [f(x^{-1}g)]\, d\mu_{G/Z} (\dot{x})\ , \ f\in \cind{J}{G}\,
       \lambda \ , \ g\in G\ .
$$

Comme pr\'ec\'edemment, les $\CC$-espaces vectoriels 
$$
M_\VV ={\rm Hom}_J\, (\lambda ,\pi )\text{ et } {\tilde M}_\VV = {\rm
  Hom}_G\, (\cind{J}{G}\, \lambda ,\pi )
$$
\noi sont canoniquement isomorphes via
$$
\Psi~: {\tilde M}_\VV \ni {\tilde \varphi}\mapsto {\tilde
  \varphi}\circ \alpha\in M_\VV
$$
\noi o\^u $\alpha$ est le plongement naturel $J$-\'equivariant de
$\lambda$ dans $\cind{J}{G}\, \lambda$, d\'efini comme dans le cas $J$
compact. 

\begin{lemma} L'inverse $\Phi$ de $\Psi$ est donn\'e par 
$\Phi (\varphi )={\tilde \varphi}$, o\^u 
$$
{\tilde \varphi}(f)=\frac{1}{\mu_{G/Z}(J/Z)}\int_{G/Z} \pi (x).\varphi
(f(x^{-1}))\, d\mu_{G/Z}(\dot{x})\ , \ \varphi\in M_\VV\ , \ f\in
\cind{J}{G}\, \lambda \ .
$$
\end{lemma}

\noi {\it D\'emonstration}. Laiss\'ee au lecteur.
\bigskip

Comme dans le cas o\^u $J$ est compact, les deux $\CC$-espaces
vectoriels $M_\VV$ et ${\tilde M}_\VV$ sont naturellement munis de
structures de $\HH (G,\check{\lambda})$-modules \`a droite. 

\begin{lemma} L'action naturelle de $\HH (G,\check{\lambda})$ sur
  $M_\VV$ est donn\'ee par
$$
\varphi .f =\int_{G/Z}\pi (x)\circ \varphi \circ f(x^{-1})\,
d{\mu_{G/Z}}(\dot{x})\ , \ \varphi \in {\rm Hom}_J (\lambda , \pi ) ,
\ f\in \HH (G,\check{\lambda})\ .
$$
\end{lemma}

\noi{\it D\'emonstration}. Laiss\'ee au lecteur.

\section{Isomorphismes d'alg\`ebre de Hecke et repr\'esentations temp\'er\'ees}

 L'objet de cette annexe est de d\'emontrer le Th\'eor\`eme (6.3).
  Nous allons en fait  d\'emontrer un r\'esultat plus g\'en\'eral qui
  est virtuellement d\'ej\`a contenu dans [BHK]. Nous utiliserons
  largement les concepts et notations de cet article. 

 \medskip

 Soit $\GG$ un groupe r\'eductif connexe d\'efini sur $F$ de groupe de points
 rationnels $G$. Fixons une mesure de Haar $\mu$ sur $G$. Soit $\hat
 G$ le dual unitaire de $G$, form\'e des classes d'isomorphie de
 repr\'esentations unitaires continues irr\'eductible de $G$ dans des
 espaces de Hilbert.
 Nous identifierons $\hat G$ de fa\c con canonique avec l'ensemble des
 classes d'isomorphie de repr\'esentations lisses,  irr\'eductibles,
 unitarisables,     de $G$ ([BHK] (2.13)). L'espace topologique $\hat
 G$ est muni de la mesure de Plancherel $\hat \mu$, duale de $\mu$. On
 note ${}_r {\hat G}$ le dual r\'eduit de $G$, c'est-\`a-dire le support
 de la mesure $\hat \mu$. Rappelons le r\'esulat fondamental suivant :

\begin{theorem}([Be], [Wa]) Le  support de la mesure de Plancherel  $\hat \mu$
 consiste en les classes d'isomorphie de repr\'esentations irr\'eductibles
 temp\'er\'ees.
\end{theorem}

 Pour toute paire $(K,\rho )$ form\'ee d'un sous-groupe ouvert compact
 $K$ et d'une repr\'esentation lisse irr\'eductible  $\rho$  de $K$, on
 note $\HH (G,\rho )$ l'alg\`ebre de Hecke sph\'erique attach\'ee \`a $\rho$
 et \`a la mesure de Haar $\mu$, comme dans [BK]{\S}4. Elle est munie
 d'une structure canonique d'alg\`ebre de Hilbert, de ${\mathcal
   C}^*$-star alg\`ebre r\'eduite  associ\'e not\'ee 
${}_r {\mathcal C}^* (G, \rho  )$ (voir [BHK] {\S}3 pour la notion
 d'alg\`ebre de Hilbert et [BK](4.3), page 152--156, pour la structure
 naturelle d'alg\`ebre de Hilbert sur une alg\`ebre de Hecke
 sph\'erique).
   Cette derni\`ere est
 obtenue de $\HH (G,\rho )$ par un  proc\'ed\'e de compl\'etion.  
 Le dual de ${}_r {\mathcal C}^* (G, \rho  )$  est not\'e 
${}_r {\hat {\mathcal C}}^*  (G,\rho )$. L'application de ``restriction
 des scalaires'' 
$$
{}_r {\hat {\mathcal C}}^*  (G,\rho ) \lra \HH (G,\rho )-\text{Mod}
$$
\noi est injective et permet d'identifier ${}_r {\hat {\mathcal C}}^*
(G,\rho )$  \`a un sous-ensemble de l'ensemble des classes
 d'isomorphie de $\HH (G,\rho )$-modules simples.

  Pour toute  repr\'esentation lisse $(\pi ,V)$ de $G$, on note 
$m_\rho (V) = V_\rho ={\rm Hom}_K\, (\rho ,V)$. On note ${}_r {\hat G}(\rho )$
 l'ensemble des $(\pi ,V)\in {}_r {\hat G}$ tels que $V_\rho \not=
 0$.  On a alors le r\'esultat suivant.

\begin{theorem} (Th\'eor\`eme B de [BHK]).  Pour tout $(\pi ,V)\in
  {}_r {\hat G}(\rho )$, la structure de  $\HH (G,\rho )$-module
   de  $V_\rho$ s'\'etend en une structure ${}_r {\mathcal C}(G,\rho
   )$-module et le module obtenu est simple.     L'application $(\pi
  ,V)\mapsto V_\rho $ induit un hom\'eomorphisme :
$$
{\hat m}_\rho~: \ {}_r {\hat G}(\rho )\lra {}_r
 {\hat {\mathcal C}}^*  (G,\rho )\ .
$$
\end{theorem}

 Donnons-nous maintenant, pour $i=1,2$, un $F$-groupe r\'eductif connexe
 ${\mathbb G}_i$, de groupe de points rationnels $G_i$, une mesure de
 Haar $\mu_i$, et deux paires $(K_i ,\rho_i )$ comme
 ci-dessus. Supposons donn\'e un isomorphisme de $\CC$-alg\`ebres
$$
j~: \ \HH (G_1 ,\rho_1 )\lra \HH (G_2 ,\rho_2 )
$$
\noi compatible avec les structures d'alg\`ebres de Hilbert. Alors $j$
s'\'etend en un isomorphisme de ${\mathcal C}^*$-star alg\`ebres
$$
{}_r {\mathcal C}^*  (G_1 ,\rho_1 )\lra {}_r  {\mathcal
    C}^*  (G_2 ,\rho_2 )
$$
\noi encore not\'e $j$,  et induit un hom\'eomorphisme

$$\begin{array}{ccc}
{\hat j}~: \ {}_r {\hat {\mathcal C}}^*  (G_1 ,\rho_1 ) & \lra & {}_r
{\hat {\mathcal
    C}}^*  (G_2 ,\rho_2 )\\
(\pi ,H )& \mapsto & (\pi\circ j ,H)
\end{array}
$$
\noi En particulier l'application naturelle
$$
j^*~: \ \HH (G,\rho_1 )-\text{Mod}\lra \HH (G,\rho_2 )-\text{Mod}
$$
\noi se ``restreint'' en ${\hat j}$.

 Supposons \`a pr\'esent que, pour $i=1,2$, $(K_i,\rho_i )$ soit un
 type au sens de Bushnell et Kutzko [BK2]. Pour $i=1,2$, notons
 ${\mathcal R}_{\rho_i}(G_i )$ la cat\'egorie des repr\'esentations lisses
 de $G_i$ engendr\'ees par leur composante $\rho_i$-isotypique. Soit 
$$
{\mathcal E}_j~:\ {\mathcal R}_{\rho_1}(G_1 )\lra   {\mathcal R}_{\rho_2}(G_2 )
$$
\noi l'\'equivalence de cat\'egories telle que le diagramme suivant
commute :
$$
\begin{array}{ccccc}
 & & {\mathcal E}_j & & \\
 & {\mathcal R}_{\rho_1}(G_1 )& \lra  &  {\mathcal R}_{\rho_2}(G_2 )& \\
m_{\rho_1} & \downarrow & & \downarrow & m_{\rho_2} \\
 & \HH (G,\rho_1 )-{\rm Mod} & \lra  & \HH (G,\rho_2 )-{\rm Mod} & \\
 &                        &  j^* &  &
\end{array}
$$
\noi Alors ce dernier diagramme se restreint en
$$
\begin{array}{ccccc}
 & & {\mathcal E}_j & & \\
 & {}_r \hat{G}_1 (\rho_1 ) & \lra  &  {}_r \hat{G}_2 (\rho_2 ) & \\
{\hat m} _{\rho_1} & \downarrow & & \downarrow & {\hat m}_{\rho_2} \\
 & {}_r {\hat C}(G_1 ,\rho_1 ) & \lra  & {}_r {\hat C}(G_1 ,\rho_1 ) & \\
 &                        &  \hat{j} &  &
\end{array}
$$

 En particulier si $(\pi ,V)$ est une repr\'esentation, lisse
 irr\'eductible,  unitaire, de $G_1$, alors $\pi$ est temp\'er\'ee
 si, et seulement, si ${\mathcal E}_j (\pi )$ est temp\'er\'e. Pour
 r\'esumer, on a le r\'esultat suivant :

\begin{theorem} Pour $i=1,2$, soit $(K_i ,\rho_i)$ un type de
  $G_i$. Soit
 $$
j~: \HH (G,\rho_1 )\lra \HH (G,\rho_2 )
$$
\noi  un isomorphisme d'alg\`ebres de Hilbert. Alors ${\mathcal E}_j$ envoie
  une repr\'esentation irr\'eductible temp\'er\'e de ${\mathcal
    R}_{\rho_1}(G_1)$ sur une repr\'esentation irr\'eductible
  temp\'er\'ee de ${\mathcal R}_{\rho_2}(G_2 )$. 
\end{theorem}

 Appliquons \`a pr\'esent ce r\'esultat au contexte du Th\'eor\`eme
 (6.3). Ici, on a:
\medskip

 -- $G_1 ={\rm Res}_{K/F}\, {\rm GL}(e,K)$ et $(K_1 ,\rho_1
 )=(I,{\mathbf 1}_I )$,
\smallskip

 -- $G_2 ={\rm GL}(N,F)$ et $(K_2 ,\rho_2 )= (J' ,\lambda ')$
\smallskip

 -- $j$ est l'isomorphisme $\Psi$,  qui est suppos\'e {\it unitaire}.
\medskip

 Or par le corollaire [BK] (5.6.7), un isomorphisme unitaire
 d'alg\`ebres de Hecke sph\'eriques  est compatible avec les
 structures unitaires de ces alg\`ebres. De plus,  par [BK] (5.6.19), page
 194, un tel isomorphisme est compatible avec les involutions. Il
 s'ensuit qu'un isomorphisme unitaire d'alg\`ebre de Hecke sph\'eriques
 est compatible avec les structure d'alg\`ebres de Hilbert, ce qui
 d\'emontre le Th\'eor\`eme (6.3).

 \begin{center}
 Universit\'e de Poitiers \\
       Laboratoire de Mathématiques et Applications\\
       UMR 7348 du CNRS\\
      SP2MI - T\'el\'eport 2\\
Bd M. et P. Curie BP 30179 \\
86962 Futuroscope Chasseneuil Cedex\\
France\\ 
\medskip
      
paul.broussous@math.univ-poitiers.fr
   \end{center}

\end{document}